\catcode`\^^Z=9
\catcode`\^^M=10
\output={\if N\header\headline={\hfill}\fi
\plainoutput\global\let\header=Y}
\magnification\magstep1
\tolerance = 500
\hsize=14.4true cm
\vsize=22.5true cm
\parindent=6true mm\overfullrule=2pt
\newcount\kapnum \kapnum=0
\newcount\parnum \parnum=0
\newcount\procnum \procnum=0
\newcount\nicknum \nicknum=1
\font\ninett=cmtt9

\font\ninebf=cmbx9

\font\sixbf=cmbx6
\font\ninesl=cmsl9

\font\nineit=cmti9

\font\ninerm=cmr9

\font\sixrm=cmr6
\font\ninei=cmmi9
\font\eighti=cmmi8
\font\sixi=cmmi6
\skewchar\ninei='177 \skewchar\eighti='177 \skewchar\sixi='177
\font\ninesy=cmsy9
\font\eightsy=cmsy8
\font\sixsy=cmsy6
\skewchar\ninesy='60 \skewchar\eightsy='60 \skewchar\sixsy='60
\font\titelfont=cmr10 scaled 1440
\font\paragratit=cmbx10 scaled 1200

\font\name=cmcsc10
\font\emph=cmbxti10

\font\tenmsbm=msbm10
\font\sevenmsbm=msbm7
%

%
\font\got=eufm10
\font\Got=eufm7
\font\teneufm=eufm10
\font\seveneufm=eufm7
\font\fiveeufm=eufm5
\newfam\eufmfam
\textfont\eufmfam=\teneufm
\scriptfont\eufmfam=\seveneufm
\scriptscriptfont\eufmfam=\fiveeufm

\font\tenmsam=msam10
\font\sevenmsam=msam7
\font\fivemsam=msam5
\newfam\msamfam
\textfont\msamfam=\tenmsam
\scriptfont\msamfam=\sevenmsam
\scriptscriptfont\msamfam=\fivemsam
\font\tenmsbm=msbm10
\font\sevenmsbm=msbm7
\font\fivemsbm=msbm5
\newfam\msbmfam
\textfont\msbmfam=\tenmsbm
\scriptfont\msbmfam=\sevenmsbm
\scriptscriptfont\msbmfam=\fivemsbm
\def\Bbb#1{{\fam\msbmfam\relax#1}}
\def\cz{{\kern0.4pt\Bbb C\kern0.7pt}
}
\def\ez{{\kern0.4pt\Bbb E\kern0.7pt}
}
\def\fz{{\kern0.4pt\Bbb F\kern0.3pt}}
\def\gz{{\kern0.4pt\Bbb Z\kern0.7pt}}
\def\hz{{\kern0.4pt\Bbb H\kern0.7pt}
}
\def\kz{{\kern0.4pt\Bbb K\kern0.7pt}
}
\def\nz{{\kern0.4pt\Bbb N\kern0.7pt}
}
\def\oz{{\kern0.4pt\Bbb O\kern0.7pt}
}
\def\rz{{\kern0.4pt\Bbb R\kern0.7pt}
}
\def\sz{{\kern0.4pt\Bbb S\kern0.7pt}
}
\def\pz{{\kern0.4pt\Bbb P\kern0.7pt}
}
\def\qz{{\kern0.4pt\Bbb Q\kern0.7pt}
}
\newskip\ttglue
\def\ninepoint{\def\rm{\fam0\ninerm}%
  \textfont0=\ninerm \scriptfont0=\sixrm \scriptscriptfont0=\fiverm
  \textfont1=\ninei \scriptfont1=\sixi \scriptscriptfont1=\fivei
  \textfont2=\ninesy \scriptfont2=\sixsy \scriptscriptfont2=\fivesy
  \textfont3=\tenex \scriptfont3=\tenex \scriptscriptfont3=\tenex
  \def\it{\fam\itfam\nineit}%
  \textfont\itfam=\nineit
  \def\sl{\fam\slfam\ninesl}%
  \textfont\slfam=\ninesl
  \def\bf{\fam\bffam\ninebf}%
  \textfont\bffam=\ninebf \scriptfont\bffam=\sixbf
   \scriptscriptfont\bffam=\fivebf
  \def\tt{\fam\ttfam\ninett}%
  \textfont\ttfam=\ninett
  \tt \ttglue=.5em plus.25em minus.15em
  \normalbaselineskip=11pt
  \font\name=cmcsc9
  \let\sc=\sevenrm
  \let\big=\ninebig
  \setbox\strutbox=\hbox{\vrule height8pt depth3pt width0pt}%
  \normalbaselines\rm
  \def\sl{\it}}

\headline={\ifodd\pageno\rightheadline\else\leftheadline\fi}
\def\rightheadline{\ninepoint Paragraphen"uberschrift\hfill\folio}
\def\leftheadline{\ninepoint\folio\hfill Chapter"uberschrift}
\let\header=Y
\def\titel#1{\need 9cm \vskip 2truecm
\parnum=0\global\advance \kapnum by 1
{\baselineskip=16pt\lineskip=16pt\rightskip0pt
plus4em\spaceskip.3333em\xspaceskip.5em\pretolerance=10000\noindent
\titelfont Chapter \uppercase\expandafter{\romannumeral\kapnum}.
#1\vskip2true cm}\def\leftheadline{\ninepoint
\folio\hfill Chapter \uppercase\expandafter{\romannumeral\kapnum}.
#1}\let\header=N
}
\def\Titel#1{\need 9cm \vskip 2truecm
\global\advance \kapnum by 1
{\baselineskip=16pt\lineskip=16pt\rightskip0pt
plus4em\spaceskip.3333em\xspaceskip.5em\pretolerance=10000\noindent
\titelfont\uppercase\expandafter{\romannumeral\kapnum}.
#1\vskip2true cm}\def\leftheadline{\ninepoint
\folio\hfill\uppercase\expandafter{\romannumeral\kapnum}.
#1}\let\header=N
}
\def\need#1cm {\par\dimen0=\pagetotal\ifdim\dimen0<\vsize
\global\advance\dimen0by#1 true cm
\ifdim\dimen0>\vsize\vfil\eject\noindent\fi\fi}
\def\neupara#1{\par\penalty-2000
\procnum=0\global\advance\parnum by 1
\vskip1cm\noindent{\paragratit \the\parnum. #1}%
\def\rightheadline{\ninepoint\S\the\parnum.\ #1\hfill \folio}%
\vskip 8mm\noindent}
\def\Proclaim #1 #2\finishproclaim {\bigbreak\noindent
{\bf#1\unskip{}. }{\it#2}\medbreak\noindent}
%
\gdef\proclaim #1 #2 #3\finishproclaim {\bigbreak\noindent%
\global\advance\procnum by 1
{%
{\relax\ifodd \nicknum
\hbox to 0pt{\vrule depth 0pt height0pt width\hsize
   \quad \ninett#3\hss}\else {}\fi}%
\bf\the\parnum.\the\procnum\ #1\unskip{}. }
{\it#2}
\medbreak\noindent}
\newcount\stunde \newcount\minute \newcount\hilfsvar
\def\uhrzeit{
    \stunde=\the\time \divide \stunde by 60
    \minute=\the\time
    \hilfsvar=\stunde \multiply \hilfsvar by 60
    \advance \minute by -\hilfsvar
    \ifnum\the\stunde<10
    \ifnum\the\minute<10
    0\the\stunde:0\the\minute~Uhr
    \else
    0\the\stunde:\the\minute~Uhr
    \fi
    \else
    \ifnum\the\minute<10
    \the\stunde:0\the\minute~Uhr
    \else
    \the\stunde:\the\minute~Uhr
    \fi
    \fi
    }
 \def\calB{{\cal B}}

\def\calC{{\cal C}} 
\def\calE{{\cal E}} 
 \def\calH{{\cal H}}
\def\calI{{\cal I}} \def\calJ{{\cal J}}

 \def\calP{{\cal P}}
 \def\calR{{\cal R}}

\def\gota{\hbox{\got a}} 
\def\gotb{\hbox{\got b}}

\def\gotp{\hbox{\got p}} 
 
\def\goto{\hbox{\got o}} 
\def\gotm{\hbox{\got m}}

 \def\Gotp{\hbox{\Got p}}

\def\Aut{\mathop{\rm Aut}\nolimits}
\def\Stab{\mathop{\rm Stab}\nolimits}

\def\dim{\mathop{\rm dim}\nolimits}

\def\GL{\mathop{\rm GL}\nolimits}

\def\im{\mathop{\rm Im}\nolimits} \def\Im{\im}

\def\kernel{\mathop{\rm kernel}\nolimits}

\def\mod{\mathop{\rm mod}\nolimits}
\def\O{{\rm O}}
\def\U{{\rm U}}

\def\proj{\mathop{\rm proj}\nolimits}

\def\re{\mathop{\rm Re}\nolimits}
\def\Re{\re}

\def\SL{\mathop{\rm SL}\nolimits}
\def\SO{\mathop{\rm SO}\nolimits}
\def\SU{\mathop{\rm SU}\nolimits}

\def\Spec{\mathop{\rm Spec}\nolimits}

\def\boxit#1{
  \vbox{\hrule\hbox{\vrule\kern6pt
  \vbox{\kern8pt#1\kern8pt}\kern6pt\vrule}\hrule}}
\def\Boxit#1{
  \vbox{\hrule\hbox{\vrule\kern2pt
  \vbox{\kern2pt#1\kern2pt}\kern2pt\vrule}\hrule}}

\def\zwischen#1{\bigbreak\noindent{\bf#1\medbreak\noindent}}

\def\smallni{\smallskip\noindent }
\def\medni{\medskip\noindent }

\def\Isom{\mathop{\;{\buildrel \sim\over\longrightarrow }\;}}
\def\lo{\longrightarrow}

\def\loma{\longmapsto}
\def\betr#1{\vert#1\vert}
\def\spitz#1{\langle#1\rangle}
\def\imag{{\rm i}}
\def\pii{\pi {\rm i}}

\def\set#1{\bigl\{\,#1\,\bigr\}}
\def\mag{\hbox{\rm i}}
\def\square{\hbox{\hbox to 0pt{$\sqcup$\hss}\hbox{$\sqcap$}}}
\def\qed{\ifmmode\square\else{\unskip\nobreak\hfil
\penalty50\hskip3em\null\nobreak\hfil\square
\parfillskip=0pt\finalhyphendemerits=0\endgraf}\fi}
\def\pn{\the\parnum.\the\procnum}
\def\downmapsto{{\buildrel
        {\vbox{\hbox{\hskip.2pt$\scriptstyle-$}}}
        \over{\raise7pt\vbox{\vskip-4pt\hbox{$\textstyle\downarrow$}}}}}
\nopagenumbers
\nicknum=0  
\font\pro =cmss10
\let\header=N
\def\transpose#1{\kern1pt{^t\kern-1pt#1}}%

\def\tr{{\hbox{\rm tr}}}

\def\bull{{\hbox{\raise.5ex\hbox{\titelfont.}}}}
\def\bullet{\bull}
\def\konteis{2.1}

\def\Dreier{3.3}
\def\nurA{4.1}

\def\Gplus{5.1}
\def\GTr{5.2}

\def\SBdreic{5.6}
\def\ThBP{5.7}
\def\GewEin{5.8}

\def\CusV{6.1}
\def\MultDet{6.2}

\def\CCl{7.1}

\def\DefBs{8.1}
\def\ThreeL{8.2}

\def\allREL{9.1}
\def\MainT{9.2}
\def\MainC{9.3}
\def\Isat{10.1}

\def\HilbI{10.3}
\def\HilbJ{10.4}

\def\SamHilb{11.2}
\def\MinPr{11.3}

\def\NonZ{12.2}

\def\fgVan{12.4}
\def\RAND#1{\vskip0pt\hbox to 0mm{\hss\vtop to 0pt{%
  \raggedright\ninepoint\parindent=0pt%
  \baselineskip=1pt\hsize=2cm #1\vss}}\noindent}
\noindent
\centerline{\titelfont A three dimensional ball quotient}%
\def\leftheadline{\ninepoint\folio\hfill
Some ball quotients with a Calabi--Yau model}%
\def\rightheadline{\ninepoint Introduction\hfill \folio}%
\headline={\ifodd\pageno\rightheadline\else\leftheadline\fi}

\vskip 1.5cm
\leftline{\it \hbox to 6cm{Eberhard Freitag\hss}
Riccardo Salvati
Manni  }
  \leftline {\it  \hbox to 6cm{Mathematisches Institut\hss}
Dipartimento di Matematica, }
\leftline {\it  \hbox to 6cm{Im Neuenheimer Feld 288\hss}
Piazzale Aldo Moro, 2}
\leftline {\it  \hbox to 6cm{D69120 Heidelberg\hss}
 I-00185 Roma, Italy. }
\leftline {\tt \hbox to 6cm{freitag@mathi.uni-heidelberg.de\hss}
salvati@mat.uniroma1.it}
\vskip1cm
\centerline{\paragratit \rm  2011}%
\vskip5mm\noindent%
\let\header=N%
\def\imag{{\rm i}}%
{\paragratit Introduction}%
\medni
Modular varieties (compactified quotients of Hermitian domains by arithmetic subgroups)
provide interesting examples of algebraic varieties. Up to finitely many cases
they are  expected to be of general type, but there are particular examples
where this is not the case.
So, for example, there are finitely many elliptic modular curves which are elliptic curves.
Generalizing these results, we constructed
in a recent paper [FS]   many Siegel threefolds which admit a Calabi--Yau model.
Some of them seem to be new. All Euler numbers which we obtained have been non-negative.
This was a motivation for us to look for other kinds of examples.  A promising class are
the ball quotients which belong to the unitary group $\U(1,3)$. Its arithmetic subgroups are
called Picard modular groups. The theory of Siegel modular varieties is far developed.
In particular the classical theory of theta functions of very interesting examples.
We just mention Igusa's result that the modular variety related to the congruence
group $\Gamma_2[4,8]$ is embedded in $P^9\cz$ where the equations are given by the quartic
Riemann relations. This is a variety of general type whose desingularization admits many
holomorphic differential forms of top-degree.
There seems to be no comparable result in the case of Picard modular varieties.
In this paper we determine a very particular example of a Picard modular variety
of general type.
On its non-singular models there exist many holomorphic differential forms.
In a forthcoming paper we will show that one can construct Calabi-Yau manifolds by
considering quotients of this variety and resolving singularities.
\smallskip
We describe the particular modular group that is considered in this paper.
Let $\U(1,3)\subset\GL(4,\cz)$
be the unitary group of the hermitian form
$$\bar z_0 w_0-\bar z_1 w_1-\bar z_2 w_2-\bar z_3 w_3.$$
It acts on the ball
$$\calB_3=\{(z_1,z_2,z_3);\quad \betr{z_1}^2+\betr{z_2}^2+\betr{z_3}^2<1\}.$$
Let $\calE$ be the ring of Eisenstein numbers. We consider the arithmetic group
$$G_3=\U(1,3)\cap\GL(4,\calE)$$ and its congruence subgroups
$$G_3[a]=\kernel\bigl(G_3\lo \GL(4,\calE/a)\bigr)\qquad (a\in\calE,\ a\ne 0).$$
We are especially interested in the group $G_3[3]$. We shall determine the structure
of the ring of modular forms $A(G_3[3])$. This algebra has $25$ generators,
$15$ modular forms $B_i$ of weight one and ten modular forms $C_j$ of weight $2$.
Both will appear as Borcherds products. The zeros are located on certain
2-subballs of $\calB_3$. The forms $C_i$ are cuspidal. Their squares define
holomorphic differential forms on the non-singular models. We shall determine all
relations between these forms and we shall obtain the dimension formula for the spaces
of modular forms and the subspaces of cusp forms in all weights.
\smallskip
The subring $A(G_3[\sqrt{-3}])$ is related to the Segre cubic. This can be derived from
the paper [F] which is an extension of the paper [AF].
It has also been worked out in detail by
Kondo [Ko]. We will obtain the structure of $A(G_3[\sqrt{-3}])$ as a by-product. This
ring is generated by 6 forms $T_i$ of weight 3 that satisfy the relations
$$\sum_{i=1}^6 T_i=0,\quad \sum_{i=1}^6 T_i^3=0$$
of the Segre cubic in a standard representation. In this way the modular variety associated
to $G_3[3]$ appears as covering of the Segre cubic of degree $3^9$.
\smallskip
The proof uses the theory of Borcherds products. This theory has been established
for the orthogonal group $\O(2,n)$. We use the natural embedding of $\U(1,n)$ into
$\O(2,2n)$ to carry over this theory to the case of the unitary group.
\smallskip
The proof  uses also rather involved computer calculations. We used the computer algebra
systems {\pro MAGMA} and {\pro SINGULAR} to perform these calculations. We feel that it is useless to publish
any programs since hard- and software are changing rapidly. Instead of this we describe
in some detail
the way how the calculations can be done such that an interested reader can
control them by writing own programs.
\smallskip
We are very grateful to G.~Pfister who explained us some fine points of the {\pro SINGULAR}
computer algebra system. We have to thank D.~Allcock who explained to us
the structure of certain unitary groups. Finally we thank S.~Kondo for his preprint [Ko].
\neupara{Orthogonal modular forms}%
We recall some basic facts about automorphic forms on $\O(2,n)$.
A real quadratic space
$V$ is a finite dimensional real vecor space that
has been equipped with a non-degenerate bilinear form
$(\cdot,\cdot)$.
We denote the associated quadratic form by $q(a)=(a,a)/2$. We call it also the {\it norm.\/}
We extend the bilinear form to $V(\cz)=V\otimes_\rz\cz$ as $\cz$-bilinear form.
\smallskip
We assume that the signature of $V$ is $(2,n)$.
The zero quadric in $V(\cz)$ is the complex submanifold defined by $(z,z)=0$.
We consider the set
$$\{z\in V(\cz);\quad (z,z)=0,\quad (z,\bar z)>0\}.$$
This is an open subset of the zero quadric which has two connected
components. We choose one of them and denote it by $\tilde\calH$.
Let
$\O'(V)$ be the subgroup of index two of the orthogonal group $\O(V)$
that preserves the two connected components. It contains the reflections along
vectors $a$ with negative $q(a)<0$ and is generated by them.
\smallskip
Now, let $M\subset V$ be an even lattice (i.e.~$q(a)\in\gz$ for all $a\in M$.)
We denote the dual lattice of $M$ by $M'$. The discriminant group of $M$ is the finite
group $M'/M$. The quadratic form $q$ induces a finite quadratic form
$\bar q:M'/M\to\qz/\gz$.
\smallskip
The integral orthogonal group $\O(M)$ consists of all $g\in\O(V)$ that preserve
$M$. We use the notation $\O'(M)=\O(M)\cap\O'(V)$. The discriminant kernel is the
group
$$\Gamma_M:=\kernel(\O'(M)\lo \Aut(M'/M)).$$
We denote by $\calH$ the image of $\tilde\calH$ in the projective space
$P^n(V(\cz))$. It is a smooth subset and the group $\O'(V)$ acts on it.
Let $\Gamma\subset\O(M')$ be a subgroup of finite index, then $\calH/\Gamma$ carries
a structure as quasi-projective variety due to the theory of Baily--Borel.
\smallskip
An automorphic form of weight $k\in\gz$ with respect to $\Gamma$
and with respect to a character $v$ on $\Gamma$ is a holomorphic
function $f$ on $\tilde\calH$ with the properties
\smallni 1)
$f(\gamma z)=v(\gamma)f(z)$ for all $\gamma\in\Gamma$,
\hfill\break 2) $f(tz)=t^{-k}f(z)$ for all $t\in\cz^\bullet$,
\hfill\break 3) the form is regular at the cusps.
\smallni
We omit the definition of 3).
We denote the space of
automorphic forms by $[\Gamma,k,v]$ or simply by $[\Gamma,k]$
if the character $v$ is trivial.
\neupara{Borcherds products}%
Let $M$ be an even lattice with bilinear form
$(\cdot,\cdot)$ and associated quadratic form $q(x)=(x,x)/2$
of signature $(2,n)$.
Borcherds' space of obstructions
consists of all modular forms $f:\hz\to\cz[M'/M]$ with the transformation law
$(f_\alpha)_{\alpha\in M'/M}$
$$\leqalignno{
&f_\alpha(\tau+1)=e^{-2\pi\imag q(\alpha)}f_\alpha(\tau),&1)\cr
&f_\alpha\Bigl(-{1\over\tau}\Bigr)=-\sqrt{\tau\over\mag}^{2+n}
    {1\over\sqrt{\#M^*/M}}\sum_{\beta\in M^*/M}
   e^{2\pi\imag(\alpha,\beta)}
            f_\beta(\tau).
    &2)\cr
    &\hbox{\rm $f$ is holomorphic at the cusp
    infinity.\qquad}
    &3)}$$
Let $\alpha\in M'$ be an element of the dual lattice and $n<0$ a negative
number. The Heegner divisor $H(\alpha,n)\subset\calH$ is the
union of all
$$v^\perp\cap\calH\qquad(v^\perp\hbox{
\rm orthogonal complement of $v$ in }P(V(\cz)),$$
where $v$ runs through all elements $\alpha\in M'$ with
$$v\equiv\alpha\,\mod\, M\quad\hbox{\rm and}\quad q(v)=n.$$
We consider $H(\alpha,n)$ as a divisor on $\calH$ by attaching
multiplicity $1$ to all components.
We have $H(\alpha,n)=H(-\alpha,n)$, more precisely, this divisor depends
only on the image of $\alpha$ in $(M'/M)/\pm1$. It is invariant under the discriminant
kernel $\Gamma_M$ and its image in $\calH/\Gamma_M$ is a closed algebraic subvariety.
\smallskip
A fundamental
result of Borcherds states (see [AF], Theorem 5.2 for this version of Borcherds theorem):
\proclaim
{Theorem}
{Assume $n>2$. A finite linear combination
$$\sum_{\alpha\in (L'/L)/\pm1,\;n<0}
  C(\alpha,n)H(\alpha,n)\qquad(C(\alpha,n)\in\nz_{\ge 0})$$
is the divisor of an automorphic form of weight $k$ with respect to $\Gamma_M$
if for every {\emph cusp form} $f$ in the the space of obstructions,
$$f_\alpha(\tau)=\sum_{n\in\qz} a_\alpha(n)\exp(2\pii n\tau),$$
the relation
$$\sum_{n<0,\,\alpha\in L'/L}a_\alpha(-n)C(\alpha,n)=0$$
holds.
The weight of this modular form is
$$k=
  \sum_{n\in\qz,\,\alpha\in L'/L}b_\alpha(n)C(\alpha,-n),$$
where $b_\alpha(n)$ denotes the Fourier coefficients
of the Eisenstein series with the constant term
$$b_\alpha(0)=\cases{-1/2&if $\alpha=0$,\cr0&else.}$$
\smallni
{\bf Corollary.} An individual Heegner divisor $H(\alpha,n)$ is
the divisor of an au\-to\-mor\-phic form if and
only if $a_\alpha(n)=0$
for every cusp form in the the space of obstructions.
}
 konteis%
\finishproclaim
\neupara{A special case}%
We consider a special lattice which can be defined by means of
the ring of Eisenstein integers.
$$\calE:=\gz[\zeta]\qquad(\zeta=-1/2+\sqrt{-3}/2).$$
We consider the isomorphism
$$\calE/\sqrt{-3}\calE\Isom\fz_3=\gz/3\gz,\quad 1\loma 1.$$
We mention that $\zeta\equiv 1$ mod $\sqrt{-3}$, hence the image of
$\zeta$ in $\fz_3$ is 1.
\smallskip
The
quadratic form $\bar x x$ equips $\calE$ with a structure as even
lattice. The associated bilinear form is $2\Re\bar x y=\tr (\bar xy)$. (This is
isomorphic to
the root lattice $A_2$.) The dual lattice is
$$\calE':={1\over\sqrt{-3}}\calE.$$ The discriminant group
$$\calE'/\calE\cong\fz_3=\gz/3\gz$$
has order three.
The elements $1/\sqrt{-3}$ and $\zeta/\sqrt{-3}$ from $\calE'$ both give
$1$ in $\fz_3$.
\smallskip
We have to
consider the lattice of signature $(2,6)$
$$M:=\calE^{1,3}, \quad M'={1\over\sqrt{-3}}M, $$
whose underlying group is
$\calE^4$ which is equipped with
$$(x,y)=2\Re(\bar x_0y_0-\bar x_1y_1-\bar x_2 x_2 -\bar x_3y_3),\quad
q(x)=\bar x_0x_0-\bar
x_2x_2-\bar x_3x_3-\bar x_4x_4.$$
The discriminant group $M'/M$ has
order $3^4$, hence $(M'/M)/\pm1$ has order $41$.
\smallskip
We want to determine the cuspidal part of the space of obstructions.
The transformation law is
$$\leqalignno{
&f_\alpha(\tau+1)=e^{-2\pi\imag q(\alpha)}f_\alpha(\tau),&1)\cr
&f_\alpha\Bigl(-{1\over\tau}\Bigr)=-\tau^4
    {1\over\sqrt{\#M^*/M}}\sum_{\beta\in M^*/M}
   e^{2\pi\imag(\alpha,\beta)}
            f_\beta(\tau).
    &2)\cr
}$$
The functions $f_\alpha$ have the period 3. They admit an expansion
$$f_\alpha(\tau)=\sum_{\nu=1}^\infty a_\nu q^{\nu/3},
\qquad q^{1/3}:=e^{2\pii\tau/3}.$$
The smallest power of $q$ that can occur is $q^{1/3}$.
This is also the smallest power that occurs in $\eta^8$, where
$$\eta(\tau)=e^{2\pii\tau/24}\prod_{\nu=1}^\infty (1-e^{2\pii\tau\nu}).$$
We recall the transformation formula $\eta(-1/\tau)=\sqrt{\tau/\imag}\;\eta(\tau)$.
This shows that $f_\alpha/\eta^8$ is a modular form of weight 0 and hence constant. This
shows.
\proclaim
{Lemma}
{The space of obstructions for $M=\calE^{1,3}$ consists of all
vector valued functions of the form
$f_\alpha=C_\alpha\eta^8$, where $C_\alpha$ are constants that satisfy
$$C_\alpha\ne 0\Longrightarrow q(\alpha)\equiv 1/3\;\mod\; 1$$
and
$$C_\alpha=   -{1\over\sqrt{\#M^*/M}}\sum_{\beta\in M^*/M}
   e^{2\pi\imag(\alpha,\beta)}
            C_\beta.$$
{\bf Corollary.} If $f=\sum  a_\nu q^{\nu}$ is a cuspidal element
of the space of obstructions then
$$a_\nu\ne 0\Longrightarrow \nu\equiv 1/3\;\mod\; 1.$$
{\bf Corollary.}
An individual Heegner divisor $H(\alpha,n)$
is the divisor of an automorphic form if $n\not\equiv -1/3\;\mod\; 1$.
}
SpOb%
\finishproclaim
In Theorem \konteis\ we gave the formula
$$k=
  \sum_{n\in\qz,\,\alpha\in L'/L}b_\alpha(n)C(\alpha,-n),$$
for the weight of a Heegner divisor
$$\sum_{\alpha\in (L'/L)/\pm1,\;n<0}
  C(\alpha,n)H(\alpha,n)\qquad(C(\alpha,n)\in\gz)$$
in the case that the obstruction condition is satisfied. But we can consider
this number in all cases, and call it then the virtual weight of the Heegner
divisor. We compute the virtual weight in some cases. The Fourier coefficients of
the corresponding Eisenstein series can be computed by means of the formulae
in [BK].
\proclaim
{Lemma}
{The virtual weight of the Heegner divisors $H(\alpha,n)$
the three cases $(\alpha,n)$
$$\Bigl(\quad {1\over\sqrt {-3}}(0,1,0,0),\ -{1\over 3}\Bigl),
\quad \Bigl({1\over\sqrt {-3}}(0,1,1,0),\ -{2\over 3}\Bigl), \quad
\Bigl((0,1,0,0),\ -1\Bigl)$$
are $1$, $9$ and $30$. In the last two cases the Borcherds products exist, but
not in the first case.}
VirtW%
\finishproclaim
(There is an analogous result in the $\calE^{1,4}$ case [AF]: the virtual weights
there are $1/2$, $5$, $27$).
The reason that in the first case a Borcherds product cannot exist, can explained by means
of the
fact that the smallest weight of a non-vanishing automorphic form is the singular
weight, which is $2$ in our case.
\smallskip
We recall that the quadratic form on $M'$ induces a finite quadratic
form $\bar q:M'/M\to \qz/\gz$  and of course also a $\qz/\gz$-valued bilinear
form $\bar q(\alpha+\beta)-\bar q(\alpha)-\bar q(\beta)$. So we can talk about
orthogonal elements in $M'/M$.
\proclaim
{Proposition}
{Let  $\alpha_2,\alpha_2,\alpha_3$ be three pairwise orthogonal
elements of $M'/M$ that are images of vectors of $M'$ of norm $-1/3$.
Then there exists a Borcherds product on $\Gamma_M$ with divisor $H(\alpha_1,-1/3)+
H(\alpha_2,-1/3)+H(\alpha_3,-1/3)$. The weight of this modular form is\/ $3$.
}
Dreier%
\finishproclaim
{\it Proof.\/} The proof rests on the explicit computation of the space
of obstructions. The weight can be taken from Proposition \Dreier.
\neupara{Ball quotients}%
We recall some basic facts about ball quotients. Let $V$ be a finite
dimensional complex vector space which is equipped with a
hermitian form $\spitz{\cdot,\cdot}$ of signature $(1,n)$. This means
that there exists an isomorphism $V\cong \cz^{n+1}$ such that
$$\spitz{z,w}=\bar z_0w_0-\bar z_1w_1-\cdots-\bar z_nw_n.$$
A line (=one dimensional sub-vector-space of $V$) is called positive,
if it is  represented by an element of positive norm ($\spitz{z,z}>0$).
We denote by $\calB=\calB(V)\subset P(V)$ the set of all positive
lines. (As usual $P(V)$ denotes the projective space of $V$, i.e.\
the set of all lines in $V$.) In the above model the component $z_0$
of an element $z\in\cz^{n+1}$ with positive norm is different from 0.
Hence each positive line contains a unique $z$ with $z_0=1$. This identifies
$\calB(V)$ with the standard $n$-ball
$$\calB_n:=\set{z\in\cz^{n};\quad \betr{z_1}^2+\cdots+\betr{z_n}^2<1}.$$
The unitary group $\U(V)\cong\U(1,n)$ acts on $\calB$ as group
of biholomorhic automorphisms.
\smallskip
\noindent We need an arithmetic structure and consider therefore an lattice
$M\subset V$. First of all, this is a discrete additive subgroup
such that $V/M$ is compact (which means that $M$ generates $V$ as
real vector space). We assume that there exists a non-real multiplier
$a\in\cz$, $aM\subset M$. Then the  set of all these multipliers is
an order $\goto$ in an imaginary quadratic field $K$. This gives as a rational
structure on $V$. An element of $V$ is called {\it rational\/}
if it is contained in $\qz M$ which is actually a $K$-vector space.
An element of $P(V)$ is called rational if it can be represented by
a rational element. Then it can be represented also by an element
of $M$. An element of $P(V)$ is called isotropic if
it can be represented by an isotropic element
$\spitz{a,a}=0$. A {\it rational boundary point\/} is a rational isotropic element
of $P(V)$. The extended ball $\calB^*$ is the union of
$\calB$ with the  set of all rational boundary points.
\smallskip
The modular group $\U(M)$ with respect to the lattice $M$ is the subgroup
of $\U(M)$ which preserves $M$. More generally one admits
a subgroup of finite index $G\subset\U(M)$. By the theory
of Baily-Borel the quotient
$$X_G:=\calB^*/G$$
carries a structure in the form of a projective algebraic variety.
Hence the uncompactified quotient
$$X_G^0:=\calB/G$$
is a quasi projective variety.
The number of
classes of rational boundary points is finite.
\smallskip
We recall the notion of a modular form of integral weight $k$. For this we consider the inverse image
$\tilde\calB$ of $\calB$ in $V-\{0\}$. It consists of all $z\in V$ with
$\spitz{z,z}>0$. This is a connected open subset.
\smallskip
Let $G\subset\U(M)$ be a subgroup of finite index and let $v:G\to\cz^\bullet$
be a character.
A modular form of weight $k$ on $G$ with respect
to a $v$ is a holomorphic function
$f:\tilde\calB\to\cz$ with the properties
\vskip1mm
\item{1)} $f(tz)=t^{-k}f(z)$ for all $t\in\cz^\bullet$,
\item{2)} $f(\gamma z)=v(\gamma)f(z)$ for all $\gamma\in G$,
\item{3)} $f$ is regular at the cusps.
\smallni
We explain the meaning of 3).
We choose a non-zero isotropic vector
$\alpha\in M\otimes\qz$ and a vector $\alpha'\in M\otimes\qz$ such that
$\spitz{\alpha,\alpha'}=1$.
We consider
$$f_{\alpha,\alpha'}(\tau):=f(\tau\alpha+\imag\alpha').$$
The positivity condition $\spitz{\tau\alpha+\imag\alpha',\tau\alpha+\imag\alpha'}>0$
means that
$\tau$ varies in an upper half plane $\Im\tau>C\ge 0$,
It is easy to see that $f_{\alpha,\alpha'}$ has some period $N$.
(This follows from
the description of the stabilizer of a cusp which
will be given in Section 8.)
We obtain a Fourier expansion
$$f_{\alpha,\alpha'}(\tau)=\sum a_ne^{2\pi\imag n\tau/N}.$$
Regularity at the cusps means that $a_n$ vanishes
for $n<0$ (for all choices of $\alpha,\alpha'$.
Then we can define
$$f(\alpha):=a_0.$$
\proclaim
{Lemma}
{The value $f(\alpha)$
does not depend on the choice of $\alpha'$.
Moreover it depends only on the $G$-orbit of
$\alpha$ if the multipliers of $f$ are trivial.
Finally
$f(C\alpha)=\bar C^kf(\alpha).$}
nurA%
\finishproclaim
{\it Proof.\/} One considers more generally
$$f_{\alpha,\alpha'}(\tau,w):=
f(\tau\alpha+\imag\alpha'+w)$$
where $w$ is orthogonal to $\alpha,\alpha'$. This function admits a
Fourier expansion
$$f_{\alpha,\alpha'}(\tau,w)=\sum a_n(w)e^{2\pi\imag n\tau/N}.$$
The functions $a_n(w)$ are holomorphic on the whole vector space
(orthogonal complement of $\cz\alpha+\cz\beta$). The function $a_0(w)$ is
an abelian function and hence constant. (This follows also from
the description of the stabilizer of a cusp given in Section 8.)
\smallskip
Now we can prove the independence of the choice of $\alpha'$. Let $\spitz{\alpha,\alpha''}=1$.
Then $\imag(\alpha''-\alpha')$ is orthogonal to $\alpha$ and hence of the form
$C\alpha+w$ where $w$ is orthogonal to $\alpha,\alpha'$.
We have $f(\tau\alpha+\imag\alpha'')=f((\tau+C)\alpha+\imag\alpha'+w)$.
The limits for $\Im\tau\to\infty$ are the same.
\smallskip
To get the value $f(C\alpha)$ we choose $\alpha'/\bar C$ as complementary vector. We get
$$\eqalign{
f(C\alpha):=&\lim_{\Im\tau\to\infty}f(\tau C\alpha+\imag\alpha'/\bar C)=\cr
\bar C^k&\lim_{\Im\tau\to\infty}f(\tau\betr C^2\alpha+\imag\alpha')=
\bar C^k\lim_{\Im\tau\to\infty}f(\tau\alpha+\imag\alpha')=\bar C^kf(\alpha).\cr}$$
This finishes the proof of Lemma \nurA.
\qed
\proclaim
{Definition}
{Assume that the order $\goto$ of multipliers is a principal ideal ring.
A cusp is a primitive isotropic vector in $M$.}
DefC%
\finishproclaim
Hence the zero dimensional boundary points are the images of the cusps
in the projective space. The number of cusps over a zero dimensional boundary
point is the number of units in $\goto$.
\smallskip
So we defined the value $f(\alpha)$ of a modular form at each cusp. In the case that
$f$ has trivial multipliers it depends only on the $G$-orbit of the cusp.
The number of these $G$-orbits (cusp-classes) is finite.
\neupara{Restriction of orthogonal to unitary automorphic forms}%
We use the notation of the previous section, so $V$ is a hermitian space.
We can consider the complex vector space $V$ also as a real vector space and equip
it with the real bilinear form
$$(a,b)=\spitz{a,b}+\spitz{b,a}.$$
The signature of this form $(2,2n)$.
This gives us an embedding
$$\U(V,\spitz{\cdot,\cdot})\subset\O(V,(\cdot,\cdot)).$$
We complexify the vector space $V$. This means that we consider $V\otimes_\rz\cz$
where the complex multiplication is given by
$$\imag(v\otimes C):=v\otimes (\imag C).$$
We extend the action of $\O(V)$ by $\cz$-linearity to $V\otimes_\rz\cz$.
Now we consider the map
$$V\lo V\otimes\cz,\quad v\loma (\imag v)\otimes 1+v\otimes\imag.$$
This map is $\cz$-linear, as can be seen from the following computation:
\smallni
$\imag v\loma (-v)\otimes 1+(\imag v)\otimes\imag$,\hfill\break
$\imag((\imag v)\otimes 1+v\otimes\imag)=\imag v\otimes\imag+v\otimes(-1)$.
\smallni
The expressions on the right-hand-side are the same.
\smallskip
This $\cz$-linear map is also compatible with the actions of
$\U(V)\subset\O(V)$: The image of $g(v)$ for $g\in\U(V)$ and $v\in V$ in $V\otimes_\rz\cz$ is
$\imag g(v)\otimes1+g(v)\otimes\imag$. Because of $\imag g(v)=g(\imag v)$ this equals
$g(\imag v\otimes 1+v\otimes\imag)$. Recall that $g$ acts on $V\otimes_\rz\cz$ by $\cz$-linear
extension.
\smallskip

The image of $\tilde\calB$ is contained in the set
defined by $(z,z)=0$ and $(z,\bar z)>0$. We can choose the connected
component $\tilde\calH$ such that it contains the image of $\tilde\calB$.
Since the group $\U(V)$ is connected, we have
that $\U(V)$ is contained in $\SO'(V)=\SO(V)\cap \O'(V)$.
\smallskip
We can use the holomorphic embedding $\tilde\calB\hookrightarrow\tilde\calH$ to restrict
an
orthogonal modular form $F$ for some group $\Gamma\subset\O(M)$ to a holomorphic function
$f$ on $\tilde\calB$ that is given by the formula $f(z)=F(\imag z\otimes 1+z\otimes\imag)$.
The equation $F(tz)=t^{-k}F(z)$ for real $t$ gives $f(tz)=t^{-k}f(z)$
for real $t$. Since $f$ is holomorphic
this must hold also for complex $t$. This shows that $f$ is a unitary automorphic form
of the same weight with respect to the group $\Gamma\cap\U(M)$ and the restricted character.
\smallskip
We apply this to the lattice $M=\calE^{1,3}$. We use the notation
$G_3=\U(M)$.
The discrimant kernel  is the kernel of the homomorphism
$$\U(\calE^{1,3})\lo\Aut\bigl((\calE/\sqrt{-3}\calE)^4\bigr).$$
This is the $G_3[\sqrt{-3}]$ which we introduced above.
\smallskip
Basic elements of $G_3$ are reflections along vectors $b$  of
norm $\spitz{b,b}=-1$. They are defined by
$$a\loma a-(1-\eta){\spitz{b,a}\over\spitz{b,b}}b$$
where $\eta$ is a 6th root of unity (hence a power of $-\zeta)$.
They transform $b$ to $\eta b$ and act as identity on the orthogonal
complement of $b$. Their order equals the order of $\eta$. We call them
biflections, triflections or hexflections correspnding to the order of
$\eta$.
Let $G^+$ be the subgroup of $G$ generated by the hexflections.
\smallskip
For a given vector $b$ of norm -1 we always have 2 possible triflections along
$b$, since there are two primitive third roots of unity.
When we talk in the following about ``the'' triflection along $b$, we mean the
triflection with factor $\zeta$. Similarly we understand by ``the'' hexflection
along $b$ always the hexflection with the factor $-\zeta$.
\smallskip
We mention some results that can be taken from [ATC]. There they are formulated
and proved for $\calE^{1,4}$. As Allcock explained to us, the proofs there also
work for $\calE^{1,3}$. In the proof one has to replace the root lattice
E$_6$ by  A$_5$.
\proclaim
{Proposition \rm [ACT]}
{The group $G_3^+$ has index two in $G_3$. The negative of the identity is not
contained in $G_3^+$. Generators of $G_3^+$ are the $5$ hexflections along the following
vectors:
$$(0,0,1,0),\quad
(0,1,1,0),\quad
(-1,-\zeta,0,0)\quad
(0,1,0,1),\quad
(0,0,0,1).$$
}
Gplus%
\finishproclaim
We need another important result.
\proclaim
{Proposition \rm [ACT]}
{The group $G_3^+$ acts transitively on the set of primitive isotropic vectors.}
GTr%
\finishproclaim
Now we consider triflections.
Because of
$$\zeta\equiv1\; \mod\; \sqrt{-3}$$
they are elements of $G_3[\sqrt{-3}]$.
\proclaim
{Proposition \rm [ACT]}
{The group $G_3[\sqrt{-3}]$ is generated by triflections.}
TrifGen%
\finishproclaim
Finally we need the following result.
\proclaim
{Proposition \rm [ACT]}
{Let $a,b\subset\calE^{1,3}$   be two vectors with the property
$a\equiv b\;\mod\;\sqrt{-3}$. Assume that either both are primitive isotropic
or that both are of norm $-1$. Then they are equivalent mod $G_3[\sqrt{-3}]$.}
OrbWu%
\finishproclaim
These theorems have geometric applications. First we treat the full modular group
$G_3$.
\proclaim
{Proposition}
{The space $\overline{\calB_3/G_3}$ has one boundary point. The images
of orthogonal complements of integral vectors of norm $-1$ define an irreducible
divisor.}
SBdrei%
\finishproclaim
Similarly we get for the congruence group of level $\sqrt{-3}$.
\proclaim
{Proposition}
{The space $\overline{\calB_3/G_3[\sqrt{-3}]}$ has $10$ boundary
points. The images
of the orthogonal complements of integral vectors of norm $-1$ define a divisor with
$15$ components.}
SBdreic%
\finishproclaim

Let $a\in \calE^{1,3}$ be a vector of norm $-1$. The intersection of its orthogonal
complement with $\tilde\calB_3$ is the fixed point set of the triflection along the
the norm -1 vector $\sqrt{-3}a$. By a short mirror in $\tilde\calB$ we understand
the fixed point set of a triflection.
\smallskip
We observe that $\zeta a$ and $\zeta^2a$
have different orthogonal complements in $\tilde\calH_3$. But their intersections with
$\tilde\calB_3$
agree.
Now we consider in $\calE^{1,3}$ three vectors $a_1,a_2,a_3$
of norm -1 which are orthogonal in the sense that $\tr(\bar a_i,a_j)=0\equiv 3$ mod 3
for $i\ne j$. We restrict the corresponding Borcherds product (see Proposition \Dreier)
to the ball $\calB_3$.
What we have seen is
that all multiplicities of the zeros are three. Hence there exists a holomorphic
cube root of the restriction. This gives the following result.
\proclaim
{Theorem}
{Let $a_1,a_2,a_3$ be integral vectors of norm $-1$ with the property
$\tr(\bar a_i,a_j)\equiv 0\; \mod\; 3$ for $i\ne j$.
There exists a modular form of weight $1$ on the group $G_3[\sqrt{-3}]$
such that the zero divisor in $\overline{\calB_3/G[\sqrt{-3}]}$
consists of three of the $15$ divisors described in Proposition \SBdreic, namely
the images of the orthogonal complements of the $a_1,a_2,a_3$.
There are $15$ triples mod $\sqrt{-3}$. Hence we get $15$ well-defined
one-dimensional spaces of
modular forms. The group $G_3$ permutes them transitively.
}
ThBP%
\finishproclaim
As an example we can take the three pairwise orthogonal vectors
$$(0,1,0,0),\quad (0,0,1,0),\quad (0,0,0,1).$$
\proclaim
{Theorem}
{The exists a modular form of weight $1$ on the group $G_3[\sqrt{-3}]$.
The zero divisor of this form consists of three $G_3[\sqrt{-3}]$-orbits.
As representatives in the ball $\calB_3$ on can take the three divisors
$$w_1=0,\quad w_2=0,\quad w_3=0.$$
The multiplicities are one.}
GewEin%
\finishproclaim
We want to determine the cusps where the 15 forms vanish. So let
$a$ be a cusp and $a$ be an integral vector of norm -1. The image of the orthogonal
complement of $b$ in $\overline{\calB_3/G_3[\sqrt{-3}]}$ contains the cusp $a$ if
and only if there exists an $g\in G_3[\sqrt{-}]$ such that $\spitz{a,g(b)}=0$.
This condition implies $\spitz{a,b}\equiv 0$ mod $\sqrt{-3}$. The converse
is also true.
\proclaim
{Lemma}
{Assume that $a$ is a cusp and $b$ and integral vector of norm $-1$ such that
$$\spitz{a,b}\equiv 0\;\mod\; \sqrt{-3}.$$
Then the image of $a$ in $\overline{\calB_3/G_3[\sqrt{-3}]}$ is contained
in the image of the orthogonal complement of $b$.}
SpitzN%
\finishproclaim
{\it Proof.\/}
Since we can replace $a,b$ by $g(a),g(b)$ with some $g\in G_3$, we have to solve
the problem only for one distinguished $a$. After that we still may replace
$b$ by some $h(b)$ where now $h\in G_3[\sqrt{-3}]$. This reduces the problem
to finitely many cases (namely six). Each of them can be settled by hand.
\qed\smallskip
One can use this Lemma to determine the cuspidal zeros of our Borcherds
products.
\proclaim
{Proposition}
{Each of the $15$ Borcherds products (Theorem \ThBP)
vanishes at $6$ of the $10$ boundary points of
$\overline{\calB_3/G_3[\sqrt{-3}]}$.
}
Nsix%
\finishproclaim
So we far we do not have any information about the multiplier systems.
We will determine them in the next section.
\neupara{The multiplier system}%
As we have seen there exists a modular form of weight 1 on the
three ball $\calB_3$ with respect to the kernel $G_3[\sqrt{-3}]$ of the homomorphism
$$\U(\calE^{1,3})\lo\Aut\bigl((\calE/\sqrt{-3}\calE)^4\bigr).$$
We have to determine the multiplier system
Since $G_3[\sqrt{-3}]$ is generated
by the triflections along vectors $a$ of norm $-1$, we need the multipliers only for them.
  If $a$ maps to vector different  from the vectors
$$\pm(0,1,0,0),\quad\pm(0,0,1,0),\quad\pm (0,0,0,1)$$
then the multiplier is 1 since otherwise the form would vanish along the
corresponding mirror. Hence it is sufficient to determine the multiplier for
the 3 triflections $R_1,R_2,R_3$ along the three vectors above (now considered
in $\calE^{1,3}$).
\smallskip
Our method will be to restrict the form to a one-dimensional ball where
it can be identified. For this we consider a pair of cusps
$\alpha,\beta$ with the property $\spitz{\alpha,\beta}=1$.
We consider the sub-lattice
$$M=\calE \alpha+\calE \beta.$$
The hermitian form is
$$\spitz{a_1\alpha+b_1\beta,a_2\alpha+b_2\beta}=\bar a_1b_2+\bar a_2b_1.$$
The lattice $M$ is isometric equivalent to $\calE^{1,1}$.
Using the basis $\alpha,\beta$, one can identify $\U(M)$ with the
set of all matrices $M\in\GL(2,\calE)$ with the property
$$\bar M'\pmatrix{0&1\cr 1&0}M\equiv\pmatrix{0&1\cr 1&0}.$$
The subgroup
$$G_1[\sqrt{-3}]:=\kernel\bigl(\U(M)
\lo\Aut(M/\sqrt{-3}M).$$
can be considered as subgroup of $G_3[\sqrt{-3}]$ by acting trivial
on the orthogonal complement of $M$.
The set $\tilde\calB_1$ corresponds in this model to the set of all
$(z_1,z_2)$ with the property $\Re(\bar z_1 z_2) >0$. The matrix $M$ acts on $(z_1,z_2)$
by multiplication from the left (considering it as column):
If $F(z_0,z_1,z_2,z_3)$ is a modular form for $G_3[\sqrt{-3}]$ then
$$f(z_1,z_2):=F(z_1\alpha+z_2\beta)$$
is a modular form on $G_1[\sqrt{-3}]$ of the same weight and the restricted character.
We want to relate  $G_1[\sqrt{-3}]$
to an elliptic modular group acting on the upper half plane.
The modular group $G_1[\sqrt{-3}]$ corresponds to the group of all matrices $\gamma={a\,b\choose c\, d}\in\GL(2,\calE)$
with the properties
$$a\equiv d\equiv 1\mod\sqrt{-3},\ b\equiv c\equiv 0\mod\sqrt{-3},\quad
\tr(\bar ac)=\tr(\bar b d)=0,\ \bar ad+\bar cd=1.$$
The transformation formula for $f$ is
$$f(az_1+bz_2,cz_1+dz_2)=v(\gamma) f(z_1,z_2)\quad\hbox{and}\quad f(tz_1,tz_2)=t^{-1}f(z_1,z_2).$$
We consider
$$f_0(\tau):=f(-1/\sqrt{-3},\tau).$$
Then $\tau$ varies in the usual upper half plane and the transformation formula reads as
$$\eqalign{f_0\bigl({\alpha\tau+\beta\over \gamma\tau+\delta}\Bigr)&=
\pmatrix{a&b\cr c&d}(\gamma\tau+\delta)f_0(\tau)
\quad\hbox{where}\cr
\pmatrix{\alpha&\beta\cr\gamma&\delta}&=\pmatrix{d&-c/\sqrt{-3}\cr -b\sqrt{-3}&a}.\cr}$$
If ${\alpha\,\beta\choose\gamma\,\delta}$ is an integral matrix from the
Hecke group $\Gamma_1[3]\subset\SL(2,\gz)$
(defined by $\alpha\equiv\delta\equiv 1$ mod 3 and $\gamma\equiv0$ mod 3),
then the ${a\,b\choose c\,d}$ belongs to the
group in question. Hence we see that $f_0(\tau)$ is a modular form
for the Hecke group $\Gamma_1[3]$.
So far we have no information about the multipliers of $f_0$. We will determine
them under a certain assumption on $F$. The vector $\alpha-\bar\zeta\beta$ has norm 1.
It is orthogonal to the vector $\alpha+\zeta\beta$ which has norm -1. Hence
$\alpha-\bar\zeta\beta$ is a point in $\tilde\calB_1$ which is fixed
under the triflection along $\alpha+\zeta\beta$.
This triflection is contained
in $G_1[\sqrt{-3}]$. Its matrix is
$$\pmatrix{\bar\zeta&1-\zeta\cr 1-\zeta&-2\bar\zeta}=\bar\zeta A\quad
\hbox{where}\quad A=\pmatrix{1&\sqrt{-3}\cr\sqrt{-3}&-2}.$$
The matrix that corresponds to $A$ in the Hecke group is
$$B=\pmatrix{1& 1\cr-3&-2}.$$
Now we assume that $F$ vanishes at the fixed point $\alpha-\bar\zeta\beta$.
Then $f_0$ vanishes at the fixed point of $B$.
This is enough to identify $f_0$ up to a constant factor.
We recall that there is an Eisenstein series of weight 1 on the group
$\Gamma_1[3]$ whose $q$-expansion is given by
$$E:=-{1\over6}-{q\over
1-q}+\sum_{\nu=1}^\infty\left[
    {q^{3\nu-1}\over1-q^{3\nu-1}}-
    {q^{3\nu+1}\over1-q^{3\nu+1}}\right].$$
It has also an expression as theta series, namely $\vartheta(S,\tau)=-6G(\tau)$ where
$$\vartheta(S,\tau)=\sum_{g\in\gz^2}e^{2\pii S[g]\tau},\qquad S=\pmatrix{2&1\cr 1&2}.$$
The Eisenstein series has trivial multiplier system on $\Gamma_1[3]$. It has
no zero at the cusps but
it also vanishes at the fixed point.
We claim that this zero is of first order and that the zeros of $E$ are equivalent
even with respect to the principal congruence subgroup of level 3.
To prove this we count the zeros of the discriminant (cusp form of weight 12).
It has a zero of third order at each of the 4 cusps of the prinicipal
congruence subgroup of level 3. Hence it has 12 zeros and a form of weight 1 must
have one zero.
It follows that
$f_0/E$ is a modular form of weight 0 and hence constant.
So we have identified the modular form and we can
determine the multiplier of the triflection $A$. Recall that $A$ is the
product of multiplication by $\bar\zeta$ and of a transformation with trivial
multiplier. From the rule $f(tz_1,tz_2)=t^{-1}f(z_1,z_2)$ we can see that the multiplier
of the triflection is $\zeta$.
\proclaim
{Proposition}
{Let $F$ be a modular form on $G_3[\sqrt{-3}]$ with an arbitrary multiplier system.
Let $\alpha,\beta$ be two cusps with the property $\spitz{\alpha,\beta}=1$.
Assume that $F(\alpha-\bar\zeta\beta)=0$. Then the multiplier of the triflection
along $\alpha+\zeta\beta$ is $\zeta$. Moreover the values of $F$ at the cusps
$\alpha,\beta$ are non-zero and related by
$$F(\beta)=F(\alpha).$$}
CusV%
\finishproclaim
{\it Proof.\/}
It only remains to compare the values at the cusps $\alpha,\beta$. We can assume that
$f_0(\tau)=\vartheta(S,\tau)$. By definition of $f_0$ we have
$$F(-\alpha/\sqrt{-3}+\tau\beta)=\vartheta(S,\tau).$$
First
we compute the value of $F$ at the isotropic vector $\beta$. By definition
$$F(\beta)=\lim_{\Im\tau\to\infty}F(\tau\beta+\imag\beta'),\qquad\spitz{\beta,\beta'}=1.$$
We can rewrite this as
$$F(\beta)=\sqrt{3}^{-1}\lim_{\Im\tau\to\infty}F(\tau\beta-\beta'/\sqrt{-3}).$$
We choose $\beta'=\alpha$ to get
$$F(\beta)=\sqrt{3}^{-1}\lim_{\im\tau\to\infty}F(\beta\tau-\alpha/\sqrt{-3})=
\sqrt{3}^{-1}\lim_{\im\tau\to\infty}\vartheta(S,\tau)=\sqrt{3}^{-1}.$$
\smallskip
Next we compute $F(\alpha)$,
$$F(\alpha)=\lim_{\im\tau\to\infty}F(\tau\alpha+\imag\beta)=
\sqrt{3}\lim_{\im\tau\to\infty}F(\tau\alpha+\sqrt{-3}\beta).$$
Now we use
$$F(\alpha\tau+\sqrt{-3}\beta)=-{1\over\tau\sqrt{-3}}F\Bigl(-{\alpha\over\sqrt{-3}}+
\Bigl(-{1\over\tau}\Bigr)\beta\Bigr)=-{1\over\tau\sqrt{-3}}
\vartheta\Bigl(S,-{1\over\tau}\Bigr).$$
The theta inversion formula states
$$\vartheta\Bigl(S,-{1\over\tau}\Bigr)={\tau\over\imag}\sqrt{\det S}^{-1}\vartheta(S^{-1},\tau)
={\tau\over\sqrt{-3}}\>\vartheta(S^{-1},\tau).$$
We get $F(\alpha\tau+\sqrt{-3}\beta)={1\over3}\vartheta(S^{-1},\tau)$ and finally
$$F(\alpha)={\sqrt 3\over 3}=\sqrt{3}^{-1}=F(\beta).\eqno\square$$
\smallskip
Now we come back to the question of determining the characters of
the 15 Borcherds products of weight 1. We consider the
special case of Theorem \GewEin.
We treat the case $R_1$ (triflection along $(0,1,0,0)$).
We consider the two  vectors
$$\alpha=(-\zeta\sqrt{-3},-\bar\zeta,-\zeta,-\zeta),\quad \beta=(\sqrt{-3},-1,1,1).$$
They are isotropic and have the property $\spitz{\alpha,\beta}=1$. We have
$$\alpha+\zeta\beta=(0,1,0,0).$$
Since the Borcherds product vanishes along  the orthogonal complement of this
vector, we get $F(\alpha-\bar\zeta\beta)=0$. Now we can apply Proposition \CusV\ and obtain
that the multiplier of $R_1$ is $\zeta$.
\proclaim
{Theorem}
{Let $F$ be the modular form on $G_3[\sqrt{-3}]$ of weight one that
vanishes along the three mirrors $w_1=0$, $w_2=0$, $w_3=0$. Then the multiplier
$v(\gamma)$ of the triflection along a vector $a\in\calE^{1,3}$ of norm -1
is $\zeta$ if $a$ is congruent mod $\sqrt{-3}$ to one of the
vectors $\pm(0,1,0,0)$, $\pm(0,0,1,0)$, $\pm(0,0,0,1)$ and $1$ otherwise.}
MultDet%
\finishproclaim
The character can be described in more detail.
We consider the subgroup $G_3[3]$ of $G_3[\sqrt{-3}]$ that acts trivially
on $(\calE/3\calE)^3$.
One can check that
$$G_3[\sqrt{-3}]/G_3[3]\cong(\gz/3\gz)^{10}.$$
It also can be checked that there exists a character on this group which
has the same effect on
triflections as described in Theorem \MultDet.
Hence both characters must agree.
\proclaim
{Lemma}
{The $15$ forms $F$ described in Theorem \ThBP\ have trivial character on the  group
$G_3[3]$.}
TrivSU%
\finishproclaim
\neupara{The congruence group of level three}%
In this section we use the model of $\calE^{1,3}$ given by the
hermitian form
$$\spitz{a,b}=\bar a_1b_2+\bar a_2b_1-\bar a_3b_3-\bar a_4b_4.$$
This is equivalent to the form we used in the previous section.
The corresponding Gram matrix is
$$H=\pmatrix{0&1&0&0\cr 1&0&0&0\cr0&0&-1&0\cr0&0&0&-1}.$$
The group $G$ can be identified with the the set of all $g\in\GL(4,\calE)$ such
that $\bar g'Hg=H$.
We know (Proposition \GTr) that every cusp is of the form
$g(e)$, $g\in G$, where
$$e:=(1,0,0,0)$$
denotes our standard cusp.
We have to study the congruence group of level three, $G_3[3]$,
especially the action of this group on
the cusps. The problem is to decide when two cusps $a,b$ are in the same $G_3[3]$-orbit.
Of course then $a\equiv b\>\mod\>3$. But it will turn out that the converse is not true.
\smallskip
We have to determine the stabilizer $\Stab$ of the standard cusp
$e=(1,0,0,0)$ in $G$. There are three types of elements in the stabilizer.
\smallni
1) The transformation
$$(a_1,a_2,a_3,a_4)\loma (a_1+a_2\sqrt{-3},a_2,a_3,a_4).$$
2) The transvections
$$(a_1,a_2,a_3,a_4)\loma \bigl(a_1+(\betr{b_1}^2+\betr{b_2}^2)a_2/2+
\bar b_1a_3+\bar b_2 a_4,\;a_2,\;a_3+a_2b_1,\;a_4+a_2b_2).$$
3) Unitary transformations in the variables $a_3,a_4$. They are generated by
$$(a_3,a_4)\loma (a_4,a_3)\quad\hbox{and}\quad (-\zeta a_3,a_4).$$
It is not difficult to show that $\Stab$ is generated by the transformations
1),2),3).
\smallskip
Due to Proposition \GTr, the cusp classes of the group $G_3[3]$
are in one-to-one correspondence
with the double cosets in
$$G_3[3]\backslash G_3/\Stab.$$
The group $G_3[3]\backslash G_3$ can be considered as a subgroup
of $\GL(4,\calE/3\calE)$. Using Proposition \Gplus\ one can compute this
group and also the image of $\Stab$ in this group. We implemented the
groups in {\pro MAGMA} and obtained the following result.
\proclaim
{Proposition}
{The order of $G_3[3]\backslash G_3$ is $2^5\cdot 3^{12}\cdot5$. The order of
the image of $\Stab$ in this group is $2^3\cdot3^7$. Hence the number
of cusp classes of $G_3[3]$ is $4860=2^2\cdot3^5\cdot 5$.
}
CCl%
\finishproclaim
We recall that in the Baily--Borel compactification $\overline{\calB_3/G_3[3]}$ the cusps
have to be counted projectively. Hence we obtain the following result.
\proclaim
{Corollary of Proposition \CCl}
{The number of boundary points of $\overline{\calB_3/G_3[3]}$ is $810$.}
CCp%
\finishproclaim
A vector in $a\in (\calE/3\calE)^4$ is called primitive if it is not of the
form $a=\sqrt{-3}b$.
It is called isotropic if $\spitz{b,b}\equiv 0$ mod 3 where $b\in\calE^4$
is an inverse image of $a$. We denote the set of all such primitive isotropic
vectors by $\calP$. The group $G_3[3]\backslash G_3$ acts on $\calP$.
One can show (for example by computation) that $G_3$ acts
transitively on $\calP$. As a consequence every vector of $\calP$ is the image
of a cusp.
There is a natural equivariant map
$$G_3[3]\backslash G/\Stab\lo\calP,\quad g\loma g(e).$$
The number of elements of $\calP$ computes as $2^2\cdot3^4\cdot5$.
This shows that over each point of $\calP$ there are three cusp classes.
This can be also understood as follows. Consider the group
$$\Stab':=\{g\in G;\ g(e)\equiv e\;\mod\; 0\}.$$
Then we have
$$G_3[3]\backslash G/\Stab'\Isom\calP.$$
We have the following result.
\proclaim
{Proposition}
{The group $\Stab$ is a normal subgroup of index three of $\Stab'$. The element
$$A:=\pmatrix{
6\zeta+4&6\zeta+18&3&-5\zeta+5\cr
3\zeta-3&18\zeta+10&3\zeta+3&6\zeta+9\cr
-3\zeta-3&-9&\zeta-2&4\zeta-1\cr
-3\zeta-3&\zeta-10&\zeta-1&5\zeta\cr}$$
is contained in $\Stab'$ but not in $\Stab$.
}
Drei%
\finishproclaim
{\it Proof.\/} The first column of $A$ is mod 3 congruent to the standard
cusp $e=(1,0,0,0)$ (considered as column). Hence $A\in\Stab'$.
We prove indirectly that
$A$ is not in $\Stab$. If $A$ is in $\Stab$ then there exists an $B\in G[3]$
such that $C:=BA$ has first column $e$. Since $C$ is unitary, the first and
third column are orthogonal. This gives $c_{23}=0$.
The norm of the third column is $-1$. This gives
$$\betr{c_{33}}^2+\betr{c_{43}}^2=1.$$
Hence either  $c_{33}=0$ or $c_{43}=0$. Since $A\equiv C$ mod $3$ we get
$a_{33}\equiv 0$ mod 3 or $a_{43}\equiv 0$ mod 3. But this is obviously false.\qed
\zwischen{Computational Access to Cusp Classes}%
We have a program that selects for any two cusps $a,b$ an element $g\in G_3$ with
$b=g(a)$.
\smallskip
Hence we can find for any cusp $a$ the corresponding double coset in $G[3]\backslash G/\Stab$.
One can compute a system of representatives. In this way one can construct an explicit system $\calC$
of representatives of the 4860 cusp classes and we have an explicit bijective map
map
$$\calC\Isom G[3]\backslash G/\Stab.$$
This makes it possible to find for each cusp $a$ its representative in
$\calC$. One just writes $a$ in in the form $a=g(e)$ to obtain an element of $G[3]\backslash G/\Stab$
and takes its inverse image in $\calC$.
\smallskip
In this way we can describe the action of $G$ on the set of representatives $\calC$.
(This factors through $G[3]$.)
\neupara{First relations}%
As in the previous section we use the hermitian form
$$\spitz{a,b}=\bar a_1b_2+\bar a_2b_1-\bar a_2b_3-\bar a_4b_4.$$
We give a list of representatives of
$G_3[\sqrt{-3}]$-orbits of pairs $\pm a$ where
$a$ is an integral vector of norm $-1$.
They correspond to the $15$ short mirrors
in $\overline{\calB_3/G_3[\sqrt{-3}]}$.
\medni
{
\halign{\qquad$#$\quad&$#$\quad&$#$\quad&$#$\quad&$#$\quad\cr
1&2&3&4&5\cr
(0,0,1,0)&(0,0,0,1)&(1,0,1,0)&(1,0,-1,0)&(1,0,0,1)\cr
\noalign{\vskip2mm}
6&7&8&9&10\cr
(1,0,0,-1)&(0,1,1,0)&(0,1,-1,0)&(0,1,0,1)&(0,1,0,-1)\cr
\noalign{\vskip2mm}
11&12&13&14&15\cr
(\zeta,-1,1,1)&(\zeta,-1,1,-1)&(\zeta,-1,-1,1)&(\zeta,-1,-1,-1)&(\zeta,1,0,0)\cr}
\medni
We give the list of triples of short mirrors which are pairwise orthogonal in the sense of
Theorem \ThBP:
\proclaim
{Definition}
{The modular forms of weight one which vanish along the following triples
of short mirrors
\smallni
\vbox{\rm\noindent
(1,2,15),
(2,4,8),
(2,3,7),
(1,6,10),
(1,5,9),
(12,13,15),
(11,14,15),
(4,6,11),
(4,5,12),
(8,10,14),
(8,9,13),
(3,6,13),
(3,5,14),
(7,10,12),
(7,9,11)}
\smallni
are denoted by
$B_1,\dots,B_{15}$
in this ordering.}
DefBs%
\finishproclaim
The forms $B_i$ are of course only determined up to constant factors. Later we will
normalize them in a suitable way.
\smallskip
It is possible to determine the values of these 15 modular forms at the 4860 cusp classes.
Let $F$ be one of the 15 forms. First one can determine the set of
cusps at which $F$ does not vanish.
Then one can decompose this set into orbits under the group $G_3[\sqrt{-3}]$.
It turns out that there are 4 orbits. It is enough to compute the value for one element in each
orbit. For this we constructed for each pair of orbits $\O_1$ and $\O_2$ a pair of vectors
$a\in\O_1$, $b\in\O_2$ that satisfies the conditions in Proposition \CusV. Then we know
$F(a)=F(b)$. Now all values of $F$ at the cusps are determined up to a constant factor.
\smallskip
Using the values of the $B_i$ at the cusps, we can describe the action of the group
$G_3$ on them.
The group $G_3$ acts on modular forms through
$(f,\gamma)\mapsto f^\gamma$, where
$$f^\gamma(z):=f(\gamma z).$$
This is an action from the right.
Up to constant factors the functions $B_i$ are permuted under this action.
We describe the action of an element $g\in G_3$ by a list
$$\pmatrix{\sigma_1&\cdots&\sigma_{15}\cr
\varepsilon_1&\cdots&\varepsilon_{15}}.$$
This list has to be read as follows:
$$B[i]^g=\varepsilon_{\sigma_i}B[\sigma(i)].$$
So the image of $B[i]$ is described by the $i$th column of this list.
\proclaim
{Lemma}
{The forms $B_i$ can be normalized in such a way that the transformation group
corresponding to $G_3$ is generated by the following three transformations.
{\ninepoint
$$\eqalign{
&\pmatrix{13&3&12&9&5&7&10&2&1&14&15&8&4&6&11\cr
1&1&-\zeta&-\bar\zeta&\zeta&1&-1&-\bar\zeta&\bar\zeta&\bar\zeta&-
\bar\zeta&1&-\bar\zeta&-\bar\zeta&-\bar\zeta},
\cr
&\pmatrix{8&12&4&2&9&7&15&3&13&11&6&1&5&10&14\cr
-1&\bar\zeta&-\zeta&-\zeta&1&-1&\zeta&1&-1&1&-\zeta&-\bar\zeta&\bar\zeta&-1&-\bar\zeta},
\cr
&\pmatrix{1&5&4&7&6&2&3&13&9&15&11&10&14&8&12\cr
-\zeta&-1&-\zeta&-1&-1&-\bar\zeta&-\bar\zeta&-\zeta&-\zeta&-1&-1&-1&-1&-\bar\zeta&-\bar\zeta}.\cr}$$
}}
ThreeL%
\finishproclaim
Using Definition \DefBs\ and the list just before it, one can verify that,
for example, $B_1B_{13}B_{15}$ and
$B_3B_5B_{7}$ have the same zeros (9 short mirrors). Hence they agree up to
a constant factor. The constant factor can be determined since they are no cusp forms
and since we know the values at the cusps. Using the normalization described in Lemma
\ThreeL\ one gets
the following relation:
$$B_1B_{13}B_{15}=B_3B_5B_{7}.$$
If one applies Lemma \ThreeL\ one can produce the following 10 relations
\smallni
$B_1B_{13}B_{15}=B_3B_5B_{7}$,
$B_4B_{11}B_{13}=B_5B_{10}B_{12}$,
$B_1B_{12}B_{14}=B_3B_4B_{6}$,\hfill\break
$B_6B_{10}B_{15}=B_7B_{11}B_{14}$,
$B_4B_9B_{15}=B_5B_8B_{14}$,
$B_2B_{12}B_{15}=B_3B_8B_{11}$,\hfill\break
$B_1B_{9}B_{11}=B_2B_5B_{6}$,
$B_6B_8B_{13}=B_7B_9B_{12}$,
$B_2B_{13}B_{14}=B_3B_{9}B_{10}$,\hfill\break
$B_1B_8B_{10}=B_2B_4B_7$.
\medni
An essential point is that all these relations are defined over $\qz$.
This is not clear from advance since the transformations described in Lemma \ThreeL\
are not defined over $\qz$.
\smallskip
In the same way one can produce  15 quaternary relations.
In this way we can prove the following result.
\proclaim
{Proposition}
{There are $10$ relations of the form $B_iB_jB_k=B_\alpha B_\beta B_\gamma$ where the six
indices are pairwise different. One of them is $B_1B_{13}B_{15}=B_3B_5B_{7}$.
The others can be obtained by means of the action of $G_3$. There are also
15 relations of the type $B_iB_jB_kB_l=B_\alpha B_\beta B_\gamma B_\delta$
with $8$ pairwise different indices. One of them is
$B_{4}B_{6}B_{13}B_{15}=B_{5}B_{7}B_{12}B_{14}$. The others can be obtained from the action
of $G_3$. These relations are all defined over $\qz$.}
DreiVier%
\finishproclaim
There are also linear relations for the third powers $B_i^3$. They are modular
forms for  $G_3[\sqrt{-3}]$  with trivial character. The ring of these modular
forms can be determined from the paper [F]. In this paper the
ring of modular forms in the 4-dimensional case  $G_4[\sqrt{-3}]$
has been determined. Our three dimensional case occurs  as  a factor of this ring.
This can be deduced from lemma 3.7 in [F]. All what we have to know here is that
the space
of modular forms of weight $3$ with trivial character with respect to
$G_3[\sqrt{-3}]$} has dimension 5.
We also refer to the paper of Kondo [Ko] where this 5-dimensional system has been
treated in detail.
\smallskip
The values of the
$B_i^3$ at the cusps gives already a 5-dimensional space. Hence we see that  in this  space there
are no cusp forms of weight 3.
So one can verify the following relations looking at the values at the cusps.
\proclaim
{Proposition}
{The third powers of the functions $B_i$ satisfy the following ten linear relations.
\smallni
$B_{1}^3-B_{10}^3+B_{11}^3-B_{14}^3-B_{15}^3=0$,
$B_{2}^3-B_{10}^3+B_{11}^3=0$,
$B_{3}^3-B_{14}^3-B_{15}^3=0$,
$B_{4}^3-B_{10}^3-B_{14}^3=0$,
$B_{5}^3+B_{11}^3-B_{15}^3=0$,
$B_{6}^3+B_{11}^3+B_{13}^3-B_{14}^3-B_{15}^3=0$,
$B_{7}^3-B_{10}^3-B_{13}^3=0$,
$B_{8}^3-B_{10}^3-B_{13}^3+B_{15}^3=0$,
$B_{9}^3-B_{11}^3-B_{13}^3+B_{15}^3=0$,
$B_{12}^3+B_{13}^3-B_{14}^3-B_{15}^3=0$.
}
ZehnR%
\finishproclaim
We denote by
$$A(G_3[3])=\bigoplus_{k=0}^\infty [G_3[3],k]$$
the ring of modular forms with respect to the principal congruence subgroup of level 3.
Here   $[G_3[3],k]$   is the space of modular forms of weight $k$ with trivial character
and similarly
$$A(G_3[\sqrt{-3}]=\bigoplus_{k=0}^\infty [G_3[\sqrt{-3}],k].$$
We consider the subring generated by the third powers $B_i^3$.
There is a relation to the Segre cubic. Recall that it is defined as follows.
One considers in the polynomial ring $\cz[T_1,\dots,T_6]$ the ideal generated by
$T_1+\cdots+T_6$ and $T_1^3+\cdots+T_6^3$. Then the Segre cubic is
the projective threefold
associated to the graded algebra
$$\cz[T_1,\dots,T_6]/\langle T_1+\cdots+T_6,\ T_1^3+\cdots+T_6^3\rangle.$$
This algebra is normal.
\proclaim
{Proposition}
{Let
$$S:=B_1^3+\cdots+B_{15}^3.$$
The assignment
$$\eqalign{
T_1&\loma3(B_{1}^3+B_{8}^3+B_{11}^3+B_{13}^3+B_{14}^3)-S,\cr
T_2&\loma3(B_{1}^3+B_{9}^3+B_{10}^3+B_{12}^3+B_{15}^3)-S,\cr
T_3&\loma3(B_{2}^3+B_{5}^3+B_{7}^3+B_{12}^3+B_{14}^3)-S,\cr
T_4&\loma3(B_{2}^3+B_{4}^3+B_{6}^3+B_{13}^3+B_{15}^3)-S,\cr
T_5&\loma3(B_{3}^3+B_{4}^3+B_{7}^3+B_{9}^3+B_{11}^3)-S,\cr
T_6&\loma3(B_{3}^3+B_{5}^3+B_{6}^3+B_{8}^3+B_{10}^3)-S,\cr}$$
defines an isomorphism
$$\cz[T_1,\dots,T_6]/\langle T_1+\cdots+T_6,\ T_1^3+\cdots+T_6^3\rangle\Isom \cz[B_1^3,\dots,B_{15}^3].$$
Moreover, the algebra $A(G_3[\sqrt{-3}])$ is generated by the $B_i^3$.
}
SegIs%
\finishproclaim
{\it Proof.\/} This follows from the results in [F] and has also been
worked out in a self contained way by Kondo [Ko].
\smallskip
On the Segre cubic the symmetric group $S_6$ acts in an obvious way. We can see
this action also in the modular picture.
Up to constant
factors it acts on the 15 modular forms $B_i$ up to constant factors
by permutation (Lemma \ThreeL). One can compute the action of $G_3$ on the
 6 expressions $T_1,\dots,T_6$. To describe it we need a certain sign character on $G_3$.
\proclaim
{Lemma}
{The group $G_3$ has a sign character $\varepsilon:G_3\to\{1,-1\}$
that associates $-1$ to each hexflection
and to the transformation $a\mapsto -a$.}
SignC%
\finishproclaim
Now we can describe the action of $G_3$ on the $T_i$.
\proclaim
{Proposition}
{There exists a surjective homomorphism $\phi:G_3\to S_6$
with the property
$$T_i(g(z))=\varepsilon(g)T_{\phi(g)(i)}.$$
The kernel of this homomorphism is the subgroup of $G_3$
generated by $G_3[\sqrt{-3}]$ and the negative of the identity.}
phiC%
\finishproclaim
Hence we have
$$G_3/\pm (G_3[\sqrt{-3}] \cong S_6.$$
Now we obtain the following Lemma.
\proclaim
{Lemma}
{The isomorphism
$$\cz[T_1,\dots,T_6]/\langle T_1+\cdots+T_6,\ T_1^3+\cdots+T_6^3\rangle\Isom \cz[B_1^3,\dots,B_{15}^3]$$
is equivariant with respect to the homomorphism $G_3\to S_6$. Here $S_6$ acts on the variables $T_i$ by
permutation in combination with the sign character of $S_6$.}
IsoEqu%
\finishproclaim
In the remaining three sections we describe the proof of the main results.
\neupara{The algebra of modular forms}%
We consider the following 10 functions:
$$\eqalign{
C_1&=B_2B_4B_{15}/B_8,\cr
C_2&=B_2B_{13}B_{15}/B_3,\cr
C_3&=B_3B_6B_{10}/B_{14},\cr
C_4&=B_3B_5B_8/B_{15},\cr
C_5&=B_8B_{13}B_{14}/B_9,\cr
C_6&=B_5B_7B_{14}/B_{15},\cr
C_7&=B_2B_6B_{15}/B_{11},\cr
C_8&=B_1B_8B_{11}/B_2,\cr
C_9&=B_6B_{13}B_{15}/B_7,\cr
C_{10}&=B_2B_4B_6/B_1.\cr}$$
Looking at the divisors (s.~Definition \DefBs) we see that they are holomorphic. Hence they are
modular forms of weight 2.
\smallskip
We want to determine the algebraic relations between the 25 forms $B_i$, $C_j$.
For this we consider variables $X_i,Y_j$ and the natural homomorphism
$$\cz[X_1,\dots,X_{15},Y_1,\dots,Y_{10}]\lo\cz[B_1,\dots,B_{15},C_1,\dots,C_{10}].$$
The consider the action of $G_3$ on the variables $X_i$ given by the formulae described in
Lemma \ThreeL. They induce obvious transformations of the variables $Y_i$. In this way we get
an action of the group $G_3$ on the polynomial ring $\cz[X_1,\dots,Y_{10}]$ such that the
homomorphism above is $G_3$-equivariant.
\smallskip
We describe some explicit relations:
\vfill\eject\noindent
\proclaim
{Proposition}
{The following tables contains explicit relations of weight $3$, $4$, $5$ and $6$
between the forms $B_i$, $C_j$. The last entry in each line gives the order of the
$G_3$-orbit of the relation where two relations are considered to be equal if the
coincide up to a constant factor.
\smallni
{Relations of weight 3}:
\smallni
\halign{\qquad#\hfil&\quad$#$\hfil&\quad\rm#\hfil\cr
type I: & X_1X_{13}X_{15}-X_3X_5X_7&10\cr
type II: &X_2^3-X_{10}^3+X_{11}^3&10\cr
type III: &Y_2X_1-X_2X_5X_7&60\cr}
\smallni
{Relations of weight 4}:
\smallni
\halign{\qquad#\hfil&\quad$#$\hfil&\quad\rm#\hfil\cr
type I: & X_4X_6X_{13}X_{15}-X_5X_7X_{12}X_{14}&15\cr
type II: &Y_{1}X_{1}^2-X_{6}^2X_{11}X_{14}-X_{7}^2X_{10}X_{15}&90\cr
type III: &Y_1Y_2-X_2X_5X_{10}X_{15}&45\cr}
\smallni
{Relations of weight 5}:
\smallni
\halign{\qquad#\hfil&\quad$#$\hfil&\quad\rm#\hfil\cr
type I: &Y_{1}^2X_{1}-X_{6}X_{11}^2X_{14}^2-X_{7}X_{10}^2X_{15}^2&90\cr
type II: &Y_{4}X_{3}X_{9}X_{13}-Y_{7}X_{1}X_{7}X_{8}+X_{2}^2X_{4}X_{10}X_{14}&180\cr
type III: & -Y_1Y_6X_4+Y_2Y_7X_8+Y_3Y_9X_{12}&15\cr}
\smallni
{Relations of weight 6}:
\smallni
\halign{\qquad#\hfil&\quad$#$\hfil&\quad\rm#\hfil\cr
type I: &Y_4^3-X_{3}^3X_{8}^3+X_{2}^3X_{12}^3&10\cr
type II: &Y_{9}X_{5}X_{6}X_{12}X_{15}+X_{1}X_{3}^2X_{8}X_{9}^2-X_{4}X_{7}X_{13}^2X_{14}^2&90\cr
type III: &-Y_5X_4X_8X_{14}X_{15}+Y_8X_6X_7X_8X_{12}+X_5X_8^3X_9X_{13}&75\cr}
\smallni
{\bf Supplement.} If one applies an element of $G_3$ to one of these relations then
one gets up to a constant factor a relation with rational coefficients.
}
allREL%
\finishproclaim
{\it Proof.\/}
We start with the supplement. Since we know the action of $G_3$ on the $X_i$ (and as a
consequence on the $Y_i$) this can be checked using generators of $G_3$.
\smallskip
The proof of the relations uses the ring $\qz[X_1,\dots,X_{15}]$. We consider the ideal
$\calI$ that is generated by the relations which involve only the $X_i$
(relations of weight 3, type I and II and relations of weight 4 type I and their
transformed under $G_3$ multiplied by constants such that they are rational).
This ideal is rather simple and it is no problem to get a {\bf Gr\"obner} basis for it using
a computer algebra system as {\pro SINGULAR}. So one can decide whether a given polynomial
from $\qz[X_1,\dots,X_{15}]$ is contained in this ideal.
The proof of the relations -- take for example the relation $C_2B_1=B_2B_5B_7$
(weight 3, type III) -- can be given as follows. Let $\Delta=B_1\cdots B_{15}$.
It is enough to show that the modular form $\Delta(C_2B_1-B_2B_5B_7)$
vanishes. But $\Delta C_2$ can be expressed as monomial in the $B_i$ by definition
of $C_2$. Hence we have to verify a relation in the ring $\qz[X_1,\dots,X_{15}]$.
It turns out that this relation comes already from the ideal $\calI$.
In this way all the listed relations can be verified.\qed
\smallskip
Now we can formulate the main result of this paper.
\proclaim
{Theorem}
{The algebra of modular forms of $A(G_3[3])$ is generated by the forms
$B_1,\dots,B_{15}$ and $C_1,\dots,C_{10}$.
Defining relations are the $G_3$-orbits
of the relations
described in Proposition \allREL. The following dimension formula holds.
$$\dim[G_3[3],k]=\cases{
0&  for $k<0$,\cr
1& for $k=0$,\cr
10& for $k=1$,\cr
130& for $k=2$,\cr
750& for $k=3$,\cr
3115&for $k=4$,\cr
-1377+(8019/2)k-2187tk^2+(729/2)k^3&for $k>4$.\cr}$$
This algebra is defined over $\qz$, where the $\qz$-structure is generated
by these generators.
\smallskip
In weight $\geq 7$ the algebra is generated by the forms $B_i$ alone.}
MainT%
\finishproclaim
We also can determine the dimensions of the spaces of cusp forms.
\proclaim
{Theorem}
{The forms $C_i$ are cusp forms. Let $[G_3[3],k]_0$ be the space of cusp forms.
Then
$$\dim[G_3[3],k]-\dim[G_3[3],k]_0=\cases{
15& for $k=1$,\cr
120& for $k=2$,\cr
405& for $k=3$,\cr
765& for $k=4$,\cr
810& for $k>4$.}$$
}
MainC%
\finishproclaim
\neupara{The proof of the main result, Gr\"obner bases}%
The proof is a mixture of commutative algebra, computer algebra and theory of modular forms.
We denote by
$$\calJ\subset\qz[X_1,\dots,X_{15},Y_1,\dots,Y_{10}]$$
the ideal generated by the
relations described in Proposition \allREL. The ideal
$$\cz\calJ\subset\cz[X_1,\dots,X_{15},Y_1,\dots,Y_{10}]$$
is invariant under the group
$G_3$. Recall that we defined also the ideal
$$\calI\subset\qz[X_1,\dots,X_{10}]$$
that is generated by relations of weight 3 and 4. The ideal $\cz\calI$ is also
invariant under $G_3$. The ideal $\calI$ is contained in $\calJ\cap\qz[X_1,\dots,X_{10}]$
but both ideals are different. The reason is as follows. If we consider the relations
of weight 5 and type III and replace in them $Y_iY_j$ by means of the relations of weight
4 and type III, we get relations in $\qz[X_1,\dots,X_{15}]$ of weight 5 that are not
contained in $\calI$. The precise picture can be obtained with the help of computer
algebra.
{\pro SINGULAR} computes for both ideals Gr\"obner bases which enable a proof of the
following statement.
\proclaim
{Lemma}
{The ideal
$$\calI_{\hbox{\sevenrm sat}}:=\calJ\cap\qz[X_1,\dots,X_{15}]$$
has the following property. It consists of all polynomials $P$ such that
$$P\cdot X_1\cdots X_{15}\in\calI.$$
Moreover it has the following property.
Let $P\in \qz[X_1,\dots,X_{15}]$ be a polynomial such that there exists a monomial
$M=X_1^{\nu_1}\cdots X_{15}^{\nu_{15}}$ with the property $MP\in\calI_{\hbox{\sevenrm sat}}$,
then $P\in \calI_{\hbox{\sevenrm sat}}$.}
Isat%
\finishproclaim
  Gr\"obner bases also give the dimensions of the ideals.
\proclaim
{Lemma}
{The Krull dimensions of the rings
$$\qz[X_1,\dots,X_{15},Y_1,\dots,Y_{10}]/\calJ,\quad
\qz[X_1,\dots,X_{15}]/\calI,\quad \qz[X_1,\dots,X_{15}]/\calI_{\hbox{\sevenrm sat}}$$
are four.}
KrullD%
\finishproclaim
We also can get information about the Hilbert polynomials.
Let
$A=\bigoplus A_k$ be a finitely generated graded algebra
over a field $A_0$. The Hilbert series is
$$\sum_k\dim A_k t^k.$$
There exist a polynomial $H(t)$, the Hilbert polynomial, such that for a
suitable natural number $k_0$ one has
$$H(k)=\dim A_k,\qquad k\equiv 0\;\mod\; k_0,\ k>>0.$$
In the case that $A$ is generated by $A_1$, one can take $k_0=1$.
In our situation the gradings have to be defined such that the weight of $X_i$ is
1 and the weights of $Y_j$ is 2.
\proclaim
{Lemma}
{The Hilbert polynomial of $\qz[X_1,\dots,X_{15}]/\calI_{\hbox{\sevenrm sat}}$
is
$$-1377+(8019/2)k-2187k^2+(729/2)k^3.$$
The Hilbert series is
$$1+15t+120t^2+660t^3+2745t^4+8898t^5+22665t^6+44550t^7+77355t^8+\cdots.$$
In the case $k>6$ the value $H(k)$ of the Hilbert polynomial gives the correct value.
}
HilbI%
\finishproclaim
As a consequence we get also the Hilbert polynomial of the ideal $\calJ$.
\proclaim
{Lemma}
{The algebras
$$\qz[X_1,\dots,X_{15},Y_1,\dots,Y_{10}]/\calJ,\quad
 \qz[X_1,\dots,X_{15}]/\calI_{\hbox{\sevenrm sat}}$$
 agree in weight $>6$. Hence they have the same Hilbert polynomial.}
 HilbJ%
 \finishproclaim
Our next goal is to prove that $\calI_{\hbox{\sevenrm sat}}$ is a prime ideal.
Since this ideal is very involved it seems to impossible to do this by straightforward
computer algebra. Instead of this we make use of rather deep results from the
theory of modular forms.
\neupara{The proof of the main result, modular forms}%
One has to investigate the relations between the modular forms $B_i,C_j$ in more detail.
\proclaim
{Lemma}
{The algebra of modular forms
$$A(\Gamma_3[3]=\bigoplus [G_3[3],k]$$
is the normalization of the
subalgebra
$$\cz[B_1,\dots,B_{15}].$$}
NormSub%
\finishproclaim
{\it Proof.\/}
The modular forms $B_i$ have no common zero in the Baily-Borel
compactification. This can be deduced from the concrete description of their
zero divisors. It  is also an immediate consequence of the result about Segre cubic.
A general result of Hilbert implies that $A(G_3[3])$ is integral over
the image of the above homomorphism. We want to have more, namely that ist is the normalization
of the image. For this we have to show that the quotient fields agree.
The analogous statement for the group $G_3[\sqrt{-3}]$ is true. The general statement
follows be means of a Galois argument. One has to show that $G_3[\sqrt{-3}]/G_3[3]$
acts faithfully on the the image ring. This follows from the transformation formulae
above.\qed
\smallskip
The dimension formula for
space of ball-modular forms has been computed by Kato using the Selberg trace
formula [Ka]. His result gives that for $k>6$ the dimension of the space of cusp forms
as $C(k-1)(k-2)(k-3)$.
The constant $C$ can be expressed by the volume of the fundamental domain.
In our case we can obtain it as follows. We compare it with the Hilbert polynomial
of the ring $A(G_3[\sqrt{-3}])$ which can be deduced from the known structure of
this ring. We make use of the fact that the highest coefficient of the Hilbert polynomial
of $A(G_3[3])$ is just the product of the highes coefficient of the Hilbert polynomial
of $A(G_3[\sqrt{-3}]$ and the covering degree which is $3^9$. In this way one
can prove $C=729/2$. To get the dimension of the space of all
modular forms, we have to add the number of cusps $810$. So we get
$$(729/2)(k-1)(k-2)(k-3)+810.$$
This agrees with our Hilbert polynomial  $H(k)$
which we defined in Lemma \HilbI. Hence we obtain the following lemma.
\proclaim
{Lemma}
{The algebras
$$A(G_3[3]),\quad
 \qz[X_1,\dots,X_{15}]/\calI_{\hbox{\sevenrm sat}}$$
have the same Hilbert polynomial $H(k)$.}
SamHilb%
\finishproclaim
We consider the natural homomorphism
$$\qz[X_1,\dots,X_{15}]\lo A(G_3[3]).$$
Its kernel is a prime ideal $\gotp$ which contains $\calI_{\hbox{\sevenrm sat}}$.
We have the inequalities
$$\dim(\qz[X_1,\dots,X_{15}]/\calI_{\hbox{\sevenrm sat}})_k\ge
\dim(\qz[X_1,\dots,X_{15}]/\gotp)_k\ge \cz[B_1,\dots,B_{15}]_k.$$
The highest coefficient of the Hilbert polynomial of $\cz[B_1,\dots,B_{15}]$
equals the highest coefficient of the Hilbert polynomial of its normalization.
Hence the above inequalities must produce equality for the highest coefficients.
This shows that the algebras
$$\qz[X_1,\dots,X_{15}]/\calI_{\hbox{\sevenrm sat}}\quad\hbox{and}\quad
\qz[X_1,\dots,X_{15}]/\gotp$$
have the same highest coefficients. Let $\gotp_1,\dots,\gotp_n$ be the minimal
prime ideals that contain $\calI_{\hbox{\sevenrm sat}}$ and such that
$\qz[X_1,\dots,X_{15}]/{\gotp_i}$ has Krull dimension 4. The ideal $\gotp$ is one of them.
It is easy to see that the algebras
$$\qz[X_1,\dots,X_{15}]/\calI_{\hbox{\sevenrm sat}}
\quad\hbox{and}\quad
\prod_{i=1}^n\qz[X_1,\dots,X_{15}]/\gotp_i$$
have the same highest coefficient. This implies $n=1$. So we have proved the following
lemma.
\proclaim
{Lemma}
{There is only one minimal prime ideal $\gotp$
containing $\calI_{\hbox{\sevenrm sat}}$ and which has the property that
$\qz[X_1,\dots,X_{15}]/\gotp$ has dimension 4.}
MinPr%
\finishproclaim
Geometrically this means that the associated projective variety has only
one irreducible 4-dimensional component. But there might be irreducible components
of smaller dimension. We want to exclude this. For this we need some results
of commutative algebra which have been developed to get a computational access to
problems as computing the primary decomposition of a polynomial ideal. We refer
to the book [GP], especially to Chapt.~4.
Following the ideas which are developed there, we choose 4 independent variables.
\proclaim
{Lemma}
{The variables $X_1,X_3,X_5,X_6$ are independent variables for
$\calI_{\hbox{\sevenrm sat}}$ in the sense that
$$\qz[X_1,X_3,X_5,X_6]\cap\calI_{\hbox{\sevenrm sat}}=0.$$}
IndVar%
\finishproclaim
We checked this with the help of computer algebra.
We want to extend the ground field $\qz$ to
the field of rational functions
$K=\qz(X_1,X_3,X_5,X_6)$
over $\qz$ in the variables $X_1,X_3,X_5,X_6$. We consider over $K$ the polynomial ring
in the remaining 11 variables,
$$\calR=K[X_2,X_4,X_7,\dots,X_{15}].$$
Then we extend $\calI_{\hbox{\sevenrm sat}}$ to this ring. We mention some general facts:
\smallni
1) The Krull dimension of $\calR/\calI_{\hbox{\sevenrm sat}}\calR$ is zero.
\smallni
2) The minimal prime ideals of $\calR$ containing
$\calI_{\hbox{\sevenrm sat}}\calR$ are in one-to-one
correspondence with the minimal prime ideals of $\cz[X_1,\dots,X_{15}]$
containing $\calI_{\hbox{\sevenrm sat}}$ and with the additional property that
the Krull dimension of their
quotient is 4.
\smallskip
We know that there is only one prime ideal with this property (Lemma \MinPr).
Hence we see that there is only one minimal prime ideal containing
$\calI_{\hbox{\sevenrm sat}}\calR$.
This must agree with the radical of $\calI_{\hbox{\sevenrm sat}}\calR$.
So we have seen that the radical
is a prime ideal. Making again use from computer algebra, one can check that
$\calI_{\hbox{\sevenrm sat}}\calR$ agrees with its radical.
In this way we see that $\calI_{\hbox{\sevenrm sat}}\calR$ is a prime
ideal. We obtain that
$$\tilde\calI=\qz[X_1,\dots,X_{15}]\cap \calI_{\hbox{\sevenrm sat}}\calR$$
is also a prime ideal. We want to prove $\tilde\calI=\calI_{\hbox{\sevenrm sat}}$.
\smallskip
The proof again uses computer algebra. The fact is that
we can construct a Gr\"obner basis of the ideal
$\calI_{\hbox{\sevenrm sat}}\calR$  with respect
to the lexicographical ordering of the variables such that $X_i>X_j$ for $i>j$.
$$\eqalign{
&(X_{1}^3)X_{15}^6+(-X_{1}^3X_{3}^3-
X_{1}^3X_{5}^3+X_{1}^3X_{6}^3)X_{15}^3+(X_{1}^3X_{3}^3X_{5}^3-
X_{3}^3X_{5}^3X_{6}^3),\cr
&X_{14}^3+X_{15}^3+(-X_{3}^3),\cr
&X_{13}^3-X_{14}^3+(-X_{5}^3+X_{6}^3),\cr
&X_{12}^3+X_{13}^3-X_{14}^3-X_{15}^3,\cr
&X_{11}^3-X_{15}^3+(X_{5}^3),\cr
&(-X_{3}X_{5}X_{6})X_{10}+(X_{1})X_{11}X_{13}X_{14},\cr
&X_{9}^3-X_{11}^3-X_{13}^3+X_{15}^3,\cr
&(-X_{3}X_{5}X_{6})X_{8}+(X_{1})X_{9}X_{12}X_{15},\cr
&(-X_{3}X_{5})X_{7}+(X_{1})X_{13}X_{15},\cr
&(-X_{3}X_{6})X_{4}+(X_{1})X_{12}X_{14},\cr
&(-X_{5}X_{6})X_{2}+(X_{1})X_{9}X_{11}.\cr}
$$
The leading terms (first term in each line) play a fundamental role in the theory
of Gr\"obner bases. We denote by $h$ the product of the coefficients of the leading terms,
$$h=-X_1^3X_3^4X_5^4X_6^4.$$
From the general theory of Gr\"obner bases one knows ([GP], Proposition 4.3.1)
that $\tilde\calI$ consists of all polynomials $P$
such the product with a suitable power of $h$ is in $\calI_{\hbox{\sevenrm sat}}$.
From Lemma \Isat\ we obtain
that $P\in\calI_{\hbox{\sevenrm sat}}$.
So we obtained the following result.
\proclaim
{Lemma}
{The ideal $\calI_{\hbox{\sevenrm sat}}$ is a prime ideal. It is the kernel
of the natural homomorphism
$$\qz[X_1,\dots,X_{15}]\lo A(G_3[3]).$$
}
IsatP%
\finishproclaim
Our next goal is to prove that the ideal
$$\calJ\subset\qz[X_1,\dots,X_{15},Y_1,\dots,Y_{10}]$$
is a prime ideal, or equivalently, the natural homomorphism
$$\qz[X_1,\dots,X_{15},Y_1,\dots,Y_{10}]/\calJ\lo A(G_3[3])$$
is injective. For this we need a detailed investigation of the relations between the
modular forms $B_i,C_j$. We make use of the fact that monomials in them
are modular forms with respect to $G_3[\sqrt{-3}]$ and that we can determine their
characters. Since modular forms with different characters are linearly independent,
we have a method to prove that all relations com from the ideal $\calJ$.
On this way we could
prove that all relations $\le 6$ in the weights are contained in the ideal $\cal J$.
Moreover we got the following dimensions
$$15,\quad 130,\quad750,\quad3115,\quad 9558,\quad22680.$$
The last two values are exactly the values H(5) and H(6) of the Hilbert polynomial.
\smallskip
In the weight $k\ge 7$ all monomials in the $C_i$, $B_i$ can be expressed
as polynomials in the $B_i$ as consequence of the relations in weight $\le 6$.
This detailed investigation of relations implies the following result.
\proclaim
{Lemma}
{The ideal $\calJ$ (see Proposition \allREL) is a prime ideal. It is the
kernel of the natural homomorphism
$$\qz[X_1,\dots,X_{15},Y_1,\dots,Y_{15}]\lo A(G_3[3]).$$}
calJP%
\finishproclaim
\neupara{The proof of the main result, normality}%
We  have to study the algebra
$$R=\qz[B_1,\dots,B_{15},C_1,\dots,C_{10}]\otimes_\qz\cz.$$
\proclaim
{Theorem}
{The ring
$$R=\qz[B_1,\dots,B_{15},C_1,\dots,C_{10}]\otimes_\qz\cz$$
is an integral domain.
The natural homomorphism
$$\qz[B_1,\dots,B_{15},C_1,\dots,C_{10}]\otimes_\qz\cz\lo A(G_3[3])$$
induces an isomorphism in weights $k\equiv 0$ mod $k_0$ for suitable $k_0$.
(The same is true already for the algebra $\qz[B_1,\dots,B_{15}]$).}
TensIsW%
\finishproclaim
{\it Proof.\/}
We only have to show that $R$ is an integral domain.
(Then a dimension argument shows that $R\to A(G_3[3])$ is injective and the Theorem
follows form the comparison of the Hilbert polynomials, Lemma \SamHilb.)
Let $\gotp_1,\dots,\gotp_n$ be the minimal prime ideals of $R$.
Their intersection is the radical.
Since $\cz$ is a separable field extension of $\qz$ the algebra $R$ is reduced.
Hence it suffices to show that $n=1$.
They primes $\gotp_i$ are conjugate under
the automorphism group of $\cz$.
We want to compare highest coefficients of Hilbert polynomials.  The
Hilbert polynomial $H(k)$ of $R$ is the same as described in Lemma \HilbI\ and
Lemma \HilbJ.
We denote its highest coefficient by $d$  ($=729/2$). The highest coefficient
of the Hilbert polynomials of $R/\gotp_i$ are $d/n$. The kernel of the homomorphism
$R\to A(G_3[3])$ is one of the primes $\gotp_i$. Hence the highest coefficient of the
Hilbert polynomial of the image is $d/n$. The highest coefficient of the
Hilbert polynomial of the normalization doesn't change.
We know that the normalization is the full ring of modular forms.
We know that the highest coefficient of the Hilbert polynomial of
$A(G_3[3])$ is $d$.
This shows $n=1$.\qed
\smallskip
Let $R$ be a noetherian local ring. The {\it depth\/} of $R$ is largest number $n$ such there
exists a regular sequence $a_1,\dots,a_n$. This means that the $a_i$ are
elements of the maximal ideal and that $a_{i+1}$ is a non-zero divisor
in $R/(a_1,\dots,a_i)$ for all $i<n$.
Recall that a noetherian ring $R$ satisfies Serre's condition $S_m$
if all localizations
by prime ideal $R_{\Gotp}$ with the property $\dim R_{\Gotp}\ge m$ have depth $\ge m$.
Serre's normality criterion says that a noetherian ring $R$ is normal
if it satisfies $S_2$ and if $R_{\Gotp}$ is regular for all prime ideals with the property
$\dim R_{\Gotp}\le 1$.
\proclaim
{Lemma}
{Let $\gotp$ be a prime ideal of a ring $R$ and $f\in\gotp$ a non-zero
divisor in $R$. Then its image $f/1$ in $R_{\Gotp}$ is a non-zero divisor.}
NonZ%
\finishproclaim
The proof is simple and can be omitted.\qed
\smallskip
We are interested in finitely generated graded algebras $R=\bigoplus R_n$
over a field $R_0$. If $f$ is a homogenous element of $R$ than, besides the usual
localization $R_f$, one considers $R_{(f)}$ which is the subring of all homogenous
elements
of degree $0$. There is  a natural open embedding $\Spec R_{(f)}\to\proj R$.
Let $\gotp$ be a prime ideal of $R$. We assume that
$\gotp$ is different from $R_{>0}$. We denote be $\gotp_0$ the ideal  generated
by the homogenous elements of $\gotp$. This is an element of $\proj R$.
\proclaim
{Lemma}
{Let $\gotp$ be a prime ideal in $R$ which is different from $R_{>0}$.
Assume that his homogenous part $\gotp_0$ is a regular point of $\proj R$.
Then $\gotp$ is a regular point of $\Spec R$.}
RegProj%
\finishproclaim
{\it Proof.\/} By assumption there exists an element $f$ of positive degree which is
not contained in $\gotp$. Then $\Spec R_{(f)}$ is an open neighborhood of $\gotp_0$.
Since the regular locus is open we can choose $f$ such that $R_{(f)}$ is a regular
ring. Since $R_f=R_{(f)}[1/f]$ is isomorphic to the polynomial ring in one variable,
it is regular too. This shows that $R_{\Gotp}$ is regular.\qed
\smallskip
We want to apply this to the ring
$$R=\qz[B_1,\dots,B_{15},C_1,\dots,C_{10}]=\qz[X_1,\dots,X_{15},Y_1,\dots,Y_{10}]/\calJ.$$
We know that this is an integral domain.
For each cusp $s$ we can consider the homogenous ideal $\gotm$
in $R$  generated by all homogeneous elements that vanish at the cusp.
This gives a finite system $S$ of ideals.
They are prime ideals which are maximal in the system of graded ideals.
We claim that the scheme $\proj R-S$ is regular.
To prove this we can extend the ground field $\qz$ to $\cz$.  So    we obtain   the
variety $\calB_3/G_3[3]$ (without cusps) which is smooth since $G_3[3]$ acts fixed-point-free.
\smallskip
We want to apply Serre's criterion to the ring $R$. We have already seen that $R_{\Gotp}$ is
regular
for all $\gotp$ with $\dim R_{\Gotp}\le 1$. Hence we have to check the condition $S_2$.
We first check it for the ideals $\gotm\in S$.
\proclaim
{Lemma}
{Let $s$ be a cusp. There exist two homogenous elements $f,g\in R$ vanishing at the cusp $s$
such that $g$ is a non-zero divisor in $R/f$.}
fgVan%
\finishproclaim
{\it Proof.\/} We consider $f=B_1$ and $g=B_8$. There are cusps in wich both vanish.
One has to prove that $B_8$ is a non-zero divisor in $R/B_1$. To express this in a
computable ideal theoretic way we recall the definition of the ideal quotient
$$\gota:\gotb:=\{x\in R;\ x\gotb\subset\gota\}.$$
The statement above means
$$(X_1,\calJ):(X_8)=(X_1,\calJ).$$
This ideal quotient can be computed by {\pro SINGULAR}.
The same method works for all other cusps.
In this way we could prove Lemma \fgVan.
\qed\smallskip
It is worthwhile to mention that this result involves the variables $Y_i$. The analogous
statement in the ring $\qz[X_1,\dots,X_{15}]$ is false.
\smallskip
Now we can verify $S_2$ for $R$. Let $\gotp$ be a prime ideal of $R$. We can assume that
$R_{\Gotp}$ is not regular. We know then that $\gotp$ contains one of the ideals
$\gotm\in S$. (This is also true for the ideal $R_{>0}$). Using Lemma \fgVan\ in connection
with Lemma \NonZ\ we get that the depth of $R_{\Gotp}$ is $\ge 2$.
This completes the proof of the main result (Theorem \MainT) of this paper. Theorem
\MainC\ follows by an easy computation in weights $<7$, since we know the values of the
forms $B_i$ at the cusps. We recall that the forms $C_i$ are cuspidal.
\medni
{\it Final Remark.\/} The canonical weight of the 3-ball is 4.  This means
that $G$-invariant differential forms of top-degree correspond to modular forms of weight
4 (with respect to a certain character). By a general result, such a cusp form of weight
4 defines a {\it holomorphic\/} differential form on any non-singular model of $X_G$.
Hence groups $G\supset G_3[3]$ have some chance to produce Calabi--Yau manifold such that the
associated Calabi-Yau form corresponds to one of the  $C_i^2$.
We will give examples in a forthcoming
paper.
\eject
 \noindent{\paragratit References}%
\vskip0.5cm\noindent
\item{[ACT]} Allcock, A., Carlson, J., Toledo, D.:
{\it The complex hyperbolic geometry of the moduli space of cubic surfaces,\/}
J. Algebraic Geom. {\bf 11}, no. 4, 659--724 (2002)
\medni
\item{[AF]} Allcock, D.\ Freitag, E.:
{\it Cubic Surfaces and Borcherds Products,\/}
Commentarii Math. Helv. Vol. {\bf 77}, Issue 2, 270--296 (2002)
\medni
\item{[BK]} Bruinier, J., Kuss, M.:
{\it Eisenstein series attached to lattices and
modular forms on orthogonal groups,\/}
Manuscr. Math. {\bf 106}, 443-459 (2001)
\medni
\item{[F]} Freitag, E.: {\it A graded algebra related to cubic surfaces,\/}
Kyushu Journal of Math. Vol. {\bf 56}, No. 2, 299--312 (2002)
\medni
\item{[FS]} Freitag, E., Salvati--Manni, R.:
{\it On Siegel three folds with a projective Calabi--Yau model,\/}
to appear in Communications in Number Theory and Physics (2012)
\medni
\item{[GP]} Greuel, G.M. Pfister, G.:
{\it A Singular introduction to commutative algebra,\/}
Springer Verlag, Berlin, Heidelberg, New York (2002)
\medni
\item{[Ka]} Kato, S.: {\it A dimension formula for a certain space of
automorphic forms of $\SU(p,1)$,\/}
Math. Ann. {\bf 266}, 457--477 (1984)
\medni
\item{[Ko]} Kondo, S.: {\it The Segre cubic and Borcherds products,\/}
Preprint (2011)
\bye